\documentclass[11pt]{article}
\usepackage{RDGThesis20}
\usepackage[margin=2cm, top=2cm, bottom=2.5cm,a4paper]{geometry}
\usepackage{amsmath}
\usepackage[normalem]{ulem}
\usepackage{amssymb}
\usepackage{bbding}
\usepackage{authblk}
\usepackage{bbm}
\usepackage{comment} % Comments package
\usepackage{fancyhdr} %pageheaders and footers in LaTeX2e
\usepackage{fancyvrb} %read/write verbatim TeX code
\usepackage{empheq}
\usepackage{graphicx} %enhanced graphics support
\usepackage{listings} %source code printer
\usepackage{mathtools}
\usepackage[version=3]{mhchem}
\usepackage{caption, threeparttable}
\usepackage{paralist} %lists package and bullet formats
\usepackage{url} %formatting for websites and e-mail addresses
\usepackage[capitalize,nameinlink,noabbrev]{cleveref}

\usepackage{newtxtext}

\def\be{\begin{equation}}
\def\ee{\end{equation}}
\def\bes{\begin{equation*}}
	\def\ees{\end{equation*}}
\def\bea{\begin{equation} \begin{aligned}}
\def\eea{\end{aligned} \end{equation}}
\def\beas{\begin{equation*} \begin{aligned}}
		\def\eeas{\end{aligned} \end{equation*}}

\def\bi{\begin{itemize}}
	\def\ei{\end{itemize}}

\usepackage[
backend=bibtex,
doi=false,
maxnames=2,
minnames=1,
firstinits=true,
maxcitenames=2,
mincitenames=1,
url=false,
eprint=false,
isbn=false,
style=numeric,
]{biblatex}

\newtheorem{proposition}{Proposition}[section]

\newtheorem{lemma}{Lemma}[section]

\newcommand{\msout}[1]{\text{\sout{\ensuremath{#1}}}}

\begin{document}
\title{Markov matrix perturbations to optimize dynamical and entropy functionals}

\author[1]{Manuel Santos Guti\'errez \thanks{Email: \texttt{m.santosgutierrez@leicester.ac.uk}}}
%\\ email: \texttt{manuel.santos-gutierrez@weizmann.ac.il}]
\author[1]{Niccolò Zagli\thanks{Email: \texttt{nz68@leicester.ac.uk}}} 
\affil[1]{School of Computing and Mathematical Sciences, University of Leicester, Leicester, LE1 7LE, UK.}
\affil[2]{Department of Statistics, Indiana University Bloomington, Indiana, USA.}
\author[2]{Giulia Carigi\thanks{Email: \texttt{gcarigi@iu.edu}}}

\maketitle

\begin{abstract}
An important problem in applied dynamical systems is to compute the external forcing that provokes the largest response of a desired observable quantity. For this, we investigate the perturbation theory of Markov matrices in connection with linear response theory in statistical physics. We use perturbative expansions to derive linear algorithms to optimize physically relevant quantities such as: entropy, Kullback-Leibler-divergence and entropy production of Markov matrices and their related probability vectors. These optimization algorithms are applied to Markov chain representations of discrete and continuous flows in and out of equilibrium.  We consider Markov matrix representations originating from Ulam-type approximations of transfer operators and a reduced order model of a turbulent flow based on unstable periodic orbits theory. We also propose a numerical protocol to recast matrix perturbations into vector field perturbations. The results allow to physically interpret the obtained optimizing perturbations without knowledge of the underlying equations, in a data-driven way. 
\end{abstract}

%\tableofcontents

\section{Introduction}

The perturbation theory of Markov matrices has extensively been developed over the past years, drawing essential contributions from diverse scientific disciplines like probability \cite{mitrophanov2024}, population dynamics \cite{caswell2019sensitivity} and statistical physics \cite{Martinelli2004}, to name a few. The ubiquity of Markov matrices, especially in applications, can surely be associated with their simplicity in linking time-evolving systems with algebraic concepts. In fact such matrices represent the dynamics of coarse-grained dynamical systems, composed of finitely many states, and studying their linear algebraic structure provides information about their evolution and long-term properties. The dependence of these algebraic features with respect to the different model parameters is crucial to understand the sensitivity and stability of the system, and identify which states are more robust to perturbations.

A fundamental application of Markov matrices is in the numerical analysis and approximation of ergodic properties of dynamical systems \cite{ulam1964,froylandapproximating1998}. Generally, a deterministic system obeying an ordinary differential equation (ODE) is equivalent to a linear operator acting on observable function spaces, as originally realized by B.~O.~Koopman \cite{Koopman1931}. Discretizing such an operator--- or its dual counterpart, namely, the Liouville operator \cite{lasota}--- yields Markov matrix representations, whose spectral properties approximate those of the full operator. A similar procedure is carried out when dealing with stochastic differential equations (SDEs). Here the operator to approximate is the Fokker-Planck operator which describes the evolution of probability density functions on the phase space \cite{risken}. 

When an external forcing is applied to a time-evolving system its statistical properties are susceptible to changing. Linking forcing, response and the system's unforced variability is an ongoing endeavor in statistical physics \cite{marconi2008}. In this regard, linear response theory provides the fundamental fluctuation dissipation theorems applicable--- with corresponding variations--- in equilibrium, near-equilibrium \cite{kubo1966} and non-equilibrium systems \cite{gallavotti_dynamical_1995,ruelle_nonequilibrium_1998, Hairer2010, Carigi_2024}. Generally, these theorems relate the response of the forced system to perturbations with its natural fluctuations in the shape of correlation functions. However, to establish if a system exhibits linear response or not is far from being a simple feat and much effort is still put in studying different systems, both theoretically and numerically. Recently, Markov matrix discretizations of dynamical systems have been used for the computation of response operators in continuous systems out of thermodynamic equilibrium \cite{Lucarini2016,SantosJSP,lucarini2025interpretableequationfreeresponsetheory}. Novel results established that suitable discretizations of the Koopman operator lead to an efficient and interpretable application of the fluctuation dissipation theorem to investigate response properties of non-equilibrium systems \cite{zagli2025,lucarini2025generalframeworklinkingfree}. 

A natural inverse question would then be, how to perturb the system to get a desired response of the system. Originally introduced in ~\cite{Galatolo_2017,MacKay_2018}, this problem is commonly referred to as \textit{The Linear Request Problem} ~\cite{KLOECKNER_2018}. Essentially, this problem corresponds to finding optimal perturbation to manage and control statistical properties of complex chaotic systems, rather than individual trajectories. There have been recent developments in this direction \cite{antown2022optimal, froyland2023optimal, galatolo2025optimal, del2025optimal, froyland2025optimal}, starting from \cite{Antown2018} where finite-state Markov chains are considered. The authors provide a strategy to select an optimal perturbation so as to maximize the linear response of the equilibrium distribution of the system, the expectation of a specified observable, and the rate of convergence of the system to the equilibrium distribution. In this work, we aim to investigate other functionals related to the probabilistic and thermodynamic structure of dynamical systems, such as entropy, Kullback-Leibler (KL) divergence, and entropy production. These quantities encode key features of the dynamics—uncertainty, divergence from a reference process, and irreversibility—and optimizing them provides insights into the underlying mechanisms that govern non-equilibrium behavior. A first motivation for addressing this problem comes from climate science, where one would be interested in determining the forcing to the Earth system that would create more information with respect to the present climatology \cite{majda_abramov_book}. Secondly, a seemingly unrelated motivation stems from Markov Chain Monte Carlo \cite[Sec.~4.8]{pavliotisbook2014}, where one would seek to sample a prescribed probability measure by suitably modifying a Markov matrix so that its invariant measure remains unchanged, yet the rate of convergence to stationarity is maximized. Both of these problems, among other shown in the main text, are here framed in the dynamical systems perspective and are linked to Markov matrix linear algebra.

Linking perturbations at the Markov matrix level with perturbations at the vector field domain is a problem that also motivates the present work. While the path from continuous dynamics to matrices is granted by Ulam's method \cite{ulam1964,froylandapproximating1998}, the reverse process is much less discussed in the context of response to forcings. Finding optimal perturbations for statistical quantities via algebraic methods is effective, but their physical interpretation would be eased by tracing back the shape of the vector field that yields such a perturbation. Consequently, the effects of optimizing statistical quantities crucially depend on the nature of the system under scrutiny. For such reason, in this paper, low-dimensional systems are investigated to cover a range of dynamical processes: discrete maps--- Lanford map---, equilibrium and non-equilibrium stochastic dynamics--- double well potential dynamics without and with applied rotation--- and chaotic flows--- the Lorenz 63 system \cite{lorenzdeterministic1963}---. On top of this, attempting to escalate the methods to higher dimensions, coarse grained representations of turbulent flows are analyzed via transition rates between unstable periodic orbits \cite{Yalniz2021,Maiocchi2022}.

The paper is organized as follows. In Sect.~\ref{sec:pert_markov} we review the perturbation theory of Markov matrices and conceptually link it with the theory of linear response. The main theoretical results of this work are reported in Sect.~\ref{eq:optimizations}, where each of the three subsections within address the problems of optimizing the three entropy functionals mentioned earlier in this paragraph. In Sect.~\ref{sec:ulam} we briefly review how to discretize à la Ulam the Liouville (Fokker-Planck) operator of an ODE (SDE) to numerically approximate ergodic features of the system using Markov matrices. We also explain a simple algorithm of how to numerically reconstruct from perturbations of the Markov Chain an associated effective drift perturbations at the SDE level. Five systems are numerically analyzed in Sect.~\ref{sec:numerics}, all of them being of different nature. We show how the optimization problems grant physical insight on the dynamics of each system. We consider Markov representations deriving not only from an Ulam discretization of the transfer operator, but also an example of a coarse grained description of a turbulent flow provided by periodic orbit theory.  Finally, future applications of this work are discussed in Sect.~\ref{sec:conclusion}. The appendices contain technical results that complement the findings herein.

\section{Perturbation-response in Markov matrices}\label{sec:pert_markov}

% The goal of this section is to briefly review as well as give a perspective on the problem of perturbing Markov matrices. The approach taken here is that of \cite{Schweitzer1968}, where we seek to isolate the contributions of the perturbation at any order in the perturbation strength. In addition, we shall relate the formulae to dynamic systems as done in \cite{Lucarini2016,SantosJSP,Antown2018}. More recently in \cite{lucarini2025interpretableequationfreeresponsetheory}, the author relates time-modulated perturbations of Markov chains to the Green's function formalism of statistical physics.

This section develops a first-order perturbation analysis for finite-state Markov chains. We first specify the conditions under which a perturbation \(\MM_{\varepsilon}=\MM+\varepsilon\mathbf{P}\) preserves stochasticity and the transition graph, and we give sufficient conditions ensuring that the perturbed chain retains a unique invariant probability vector. 
We then compute the first-order change in the probability distribution after a finite number of steps when a constant perturbation is turned on at discrete time $n=0$ and retained at every subsequent step. By taking the long-time limit, we obtain the corresponding first-order change in the invariant distribution, expressed in terms of the fundamental matrix of the unperturbed Markov chain. The stationary perturbation formulae are closely related to the classical analysis of \cite{Schweitzer1968}, while their finite-time formulation connects with response-theoretic approaches to Markov representations of dynamical systems \cite{Lucarini2016,SantosJSP,Antown2018}. 
More general time-modulated perturbations and their Green-function representation are discussed in \cite{lucarini2025interpretableequationfreeresponsetheory}; the constant in time perturbation considered here constitutes a special case of that framework.

Consider a Markov chain on \(N\) states with transition matrix \( \MM=(M_{ij})_{i,j=1}^{N} \) with column-stochastic convention: \(M_{ij}\) is the probability of transitioning from state \(j\) to state \(i\). Hence,
\[
    M_{ij}\geq 0,
    \qquad
    \mathbf{1}^{\top}\MM=\mathbf{1}^{\top},
\]
where \(\mathbf{1}\in\mathbb{R}^{N}\) denotes the vector whose entries are all equal to one. 
Throughout this work, the unperturbed Markov matrix \(\MM\) is assumed to satisfy the conditions summarized in the following

\begin{Assumption}\label{assum:markov} Let \(\MM\) be a Markov matrix such that:
    \begin{enumerate} 
    \item[\textbf{(A1)}] \(\MM\) is \emph{primitive}; that is, there exists an integer \(m\geq 1\) such that every entry of \(\MM^m\) is strictly positive; 
    \item[\textbf{(A2)}] \(\MM\) is \emph{scrambling} in the sense of \cite[Definition~3.10]{ipsenergodicity2011}, namely, there is a positive column overlap
    \begin{equation} 
       \alpha(\MM):= \min_{1\leq i,j\leq N} \sum_{k=1}^{N} \min\{M_{ki},M_{kj}\} >0;
        \label{eq:scrambling_condition} 
    \end{equation} 
    \item[\textbf{(A3)}] No transition probability is equal to one, \( M_{ij}<1\), \(i,j=1,\ldots,N\).
    \end{enumerate} 
\end{Assumption}

    These conditions carry useful consequences. 
    First, Assumption~\textbf{(A1)} and the Perron--Frobenius theorem imply that the unit eigenvalue of \(\MM\) is algebraically simple and that \(\MM\) possesses a unique invariant probability vector \[ \uu=(u_i)_{i=1}^{N} \] whose components are strictly positive, namely
    \begin{equation} 
        \MM\uu=\uu, \qquad \mathbf{1}^{\top}\uu=1, \qquad u_i>0 \quad \text{for all }i=1,\ldots,N. \label{eq:inv_mes} 
    \end{equation} 
    Moreover, primitivity in finite-state chains, as those under analysis here, implies that \(\MM\) is mixing, namely that for every probability vector \(\xx\), 
    \begin{equation*} 
        \lim_{n\to\infty}\MM^n\xx=\uu. 
        % \label{eq:mixing_convergence} 
    \end{equation*}
   Next, Assumption~\textbf{(A2)} helps quantify the degree of mixing as it gives information on the Dobrushin's ergodicity coefficient $\tau(\mathbf{M})$. Such a coefficient is defined as 
   \begin{equation}\label{eq:erg_con}
        \tau(\MM) = \max_{ \|\xx\|_1 = 1, \mathbf{1}^{\top}\xx = 0} \|\MM\xx\|_1,
    \end{equation}
   which boils down to the matrix norm in the space of sum-zero vectors in $\mathbb{R}^{N}$. If $\MM$ is scrambling, then $\tau(\mathbf{M}) <1$. This follows using the classic equivalent formulations for $\tau$ (see e.g. \cite{ipsenergodicity2011})
    \begin{equation*}
        \tau(\MM) = \frac{1}{2}\max_{j,k} \sum_{i =1}^N |M_{ij} - M_{ik} |_1 = 1 - \alpha(\MM).
    \end{equation*}
    Moreover, if $\uu_1$ and $\uu_2$ are two probability vectors, $\uu_1 - \uu_2$ is a sum-zero vector and by \eqref{eq:erg_con}
    \begin{equation*}
        \|\MM \uu_1 - \MM\uu_2\|_1 \leq \tau(\MM) \|\uu_1 - \uu_2\|_1.
    \end{equation*}
    Therefore the scrambling property also implies that $\MM$ is a strict contraction on the probability simplex, and, iterating, it ensures exponential (or geometric) mixing
    \begin{equation}\label{eq:exp_mix}
        \|\MM^n \ww -\uu \|_1 \leq \tau(\MM)^n \|\ww - \uu\|_1 = (1 - \alpha(\MM))^n  \|\ww - \uu\|_1 \leq e^{- n \alpha(\MM)}  \|\ww - \uu\|_1,
    \end{equation}
    where $\uu$ is the invariant probability vector of $\MM$ and $\ww$ is any other initial probability vector. Another useful consequence of the scrambling property is that $\MM$ exhibits a positive spectral gap, namely the second largest eigenvalue in modulus is strictly smaller than one. In fact, if $\xx$ is an eigenvector for $\MM$ with eigenvalue $\lambda\neq 1$, then $\mathbf{1}^\top \xx = 0$ and 
    \begin{equation*}
        |\lambda| \|\xx\|_1 = \|\MM \xx\|_1 \leq \tau(\MM) \|\xx\|_1,
    \end{equation*}
    so that $|\lambda|\leq \tau(\MM)$ for all $\lambda\neq 1$.
    % giving 
    % $$
    % \rho_*(\MM):= \max_{\lambda\in \sigma(\MM), \lambda \neq 1} |\lambda| \leq \tau (\MM) <1.
    % $$
    %  or, equivalently,
    %  \begin{equation*}
    %      \operatorname{Gap}(\MM):= 1 - \rho_*(\MM) \geq 1 - \tau(\MM) = \alpha(\MM) >0. 
    %  \end{equation*}

    Observe that, for a generic vector $\xx$ and any $n\geq 0$
    \begin{align*}
        \left(\MM^n - \uu \mathbf{1}^\top\right) \xx = \MM^n \left( \xx - (\mathbf{1}^\top \xx) \uu\right)
    \end{align*}
    so that, by definition \eqref{eq:erg_con} of $\tau$, and the fact that $\uu$ is a probability vector, we have:
    \begin{equation*}
        \|\left(\MM^n - \uu \mathbf{1}^\top\right) \xx \|_1 = \|\MM^n \left( \xx - (\mathbf{1}^\top \xx) \uu\right)\|_1 \leq \tau(\MM)^n \|\xx - (\mathbf{1}^\top \xx)\uu \|_1 \leq 2  \tau(\MM)^n \|\xx\|_1.
    \end{equation*}
    This implies the following useful bound on the 1-matrix norm 
    \begin{equation}\label{eq:matrix_norm_bound}
        \|\MM^n - \uu \mathbf{1}^\top\|_1 \leq 2  \tau(\MM)^n, \quad n\geq 0.
    \end{equation}
   
   Last, Assumption~\textbf{(A3)} is independent of primitivity and scrambling. Together with column-stochasticity, it implies that every column of \(\MM\) contains at least two strictly positive entries, so \(\MM\) has at least $2N$ nonzero entries. 

We now introduce a perturbation of our Markov transition matrix and discuss the conditions it should fulfill so that the perturbed dynamics still satisfies Assumption~\ref{assum:markov}. Let \( \PP=(P_{ij})_{i,j=1}^{N}\) be a perturbation direction, and let $\eps\geq 0$ denote its amplitude, so that the perturbed transition matrix is defined by
\begin{equation}
    \MM_{\eps}:=\MM+\eps\PP.
    \label{eq:perturbed_matrix}
\end{equation}
% 
% We impose the normalization
% \begin{equation}
%     \|\PP\|_{\mathrm{F}}=1.
%     \label{eq:perturbation_normalization}
% \end{equation}
% so that \(\PP\) determines the shape and direction of the perturbation, whereas \(\eps\) determines its overall intensity:
% \[
%     \|\MM_{\eps}-\MM\|_{\mathrm{F}}
%     =
%     \eps.
% \]
The following assumptions specify the class of admissible perturbation directions considered here.
\begin{Assumption}\label{assum:perturb}
    Let \( \PP=(P_{ij})_{i,j=1}^{N}\) be a matrix such that
    \begin{itemize}
        \item[(C1)] $\|\PP\|_{F}=1$;
        \item[(C2)] $\mathbf{1}^\top \PP = \mathbf{0}^\top$ i.e. the columns of $\PP$ sum to zero 
        \item[(C3)] $ P_{ij} = 0$ if $M_{ij} = 0$ i.e. no new transitions are created as results of the perturbation.
    \end{itemize}
\end{Assumption}
The constraint (C1) ensures that the Frobenius norm is fixed, so that the bookkeeping parameter $\varepsilon$ accounts for the amplitude of the perturbation. Constraint (C2) ensures conservation of probabilities under the action of $\mathbf{P}$, as $\mathbf{1}^\top \mathbf{M}_\eps = \mathbf{1}^\top$. If, in addition to (C2), we choose $\eps$ small enough, then $\MM_\eps$ is still a Markov matrix. In fact, setting 
\begin{equation}
    \label{eq:eps_M(P)}
    0\leq \eps< \eps_\MM(\PP) = \min_{i,j\, :\, P_{ij}<0} \frac{M_{ij}}{|P_{ij}|}
\end{equation}
ensures that $\MM_\eps$ has non-negative entries.

Next we want to ensure that the perturbed matrix \(\MM_{\eps}\) also possesses a unique invariant probability vector \(\vv\) i.e.
\begin{equation}\label{eq:pert_eig}
    \MM_{\eps}\vv = (\MM+\eps\PP)\vv = \vv,
    \qquad
    \mathbf{1}^{\top}\vv=1.
\end{equation}
A sufficient condition is ensuring $\MM_{\eps}$ is also exponentially mixing, which we saw is implied by having an ergodicity coefficient smaller than one. We now show that taking the parameter $\eps$ small enough ensures the perturbation preserves this property of $\MM$. 
 % The coefficient $\tau(\mathbf{M})$, furthermore, has additivity and bounding properties--- see \cite{ipsenergodicity2011}--- that can be exploited to find a range of perturbation strengths below which the matrix $\mathbf{M}$ has a simple leading eigenvalue, without necessarily assuming diagonalizability:
\begin{proposition}\label{prop:ergodicity bound}
	Let $\MM$ be an $N\times N$ Markov matrix with $\tau(\MM)<1$. Consider the perturbation $\MM_\eps = \MM + \eps \PP$ with $\PP$ as in Assumption~\ref{assum:perturb} and $\eps$ as in \eqref{eq:eps_M(P)}. If $\eps$ satisfies also the following bound
 \begin{equation}\label{eq:bound_eps}
    % \eps < \min\left\{ \frac{\sqrt{2}(1-  \tau(\MM)}{\sqrt{N}}, \eps_\MM(\PP)\right\},
      \eps < \frac{1-  \tau(\MM)}{\sqrt{2N}},
 \end{equation}
then the ergodic coefficient of the perturbed Markov matrix is such that 
\begin{equation}
    \tau(\MM_\eps)< \frac{\tau(\MM) + 1}{2}<1
\end{equation}
for any such $\eps$. 
\end{proposition}

\begin{proof}
We start by noticing that, by the triangular inequality, $\tau$ is subadditive in its argument as 
\begin{align}
    \tau (\MM_\eps) &= \max_{ \|\xx\|_1 = 1, \mathbf{1}^{\top}\xx = 0} \|(\MM + \eps \PP)\xx\|_1 \notag\\ 
    &\leq \max_{ \|\xx\|_1 = 1, \mathbf{1}^{\top}\xx = 0} \left( \|\MM\xx \|_1  + \eps\|\PP \xx \|_1\right) = \tau(\MM) + \eps\tau(\PP).\label{eq:tau_eps:1}
\end{align}
By the definition \eqref{eq:erg_con} of $\tau$ and Assumption~(C1), it follows 
\begin{equation}\label{eq:tau_eps:2}
     \tau(\PP) \leq \|\PP\|_1\leq \left(\frac{N}{2}\right)^{1/2}\|\PP\|_F = \left(\frac{N}{2}\right)^{1/2}.
\end{equation}
Then from \eqref{eq:tau_eps:1}, \eqref{eq:tau_eps:2} and \eqref{eq:bound_eps} we have the desired result as 
\begin{equation}\label{eq:bound:tau_eps}
     \tau (\MM_\eps) \leq \tau(\MM) + \eps \left(\frac{N}{2}\right)^{1/2} < \frac{1}{2} +\frac{\tau(\mathbf{M})}{2}<1
\end{equation}
for each $\eps$ satisfying \eqref{eq:eps_M(P)}-\eqref{eq:bound_eps}.
\end{proof}

As observed for $\MM$, $\tau(\MM_\eps)<1$ implies that the perturbed dynamics is still exponentially mixing and exhibits spectral gap, both uniformly over all the $\eps$ satisfying \eqref{eq:eps_M(P)}-\eqref{eq:bound_eps}.

\begin{remark}[Diagonalizable matrices]
    In the case $\MM$ is diagonalizable, then we can ensure uniqueness of the invariant measure thanks to the Bauer-Fike theorem \cite{bauerfike} instead. This is showed for example in \cite{Mitrophanov2003}.
    % and we recall the conditions for this result here. 

    % Let $\|\cdot\|_2$ be the matrix spectral norm, and assume the Markov matrix $\MM$ can be decomposed as $\MM = \XX \Lambda \XX^{-1}$ for some diagonal matrix $\boldsymbol{\Lambda}=(\lambda_j)_{j=1}^N$ and invertible matrix $\XX$. \textcolor{red}{Without loss of generality, we order the eigenvalues in decreasing order of modulus, $\lambda_1=1>|\lambda_2|\geq\ldots\geq |\lambda_N|$.} Then the Bauer-Fike theorem guarantees that if 
    % $$\|\mathbf{X}\|_2\|\mathbf{X}^{-1}\|_2\| \PP \|_2<\frac{1}{2}\min_{2\leq j \leq N}|1 - \lambda_j|,$$ 
    % the unit eigenvalue remains isolated and simple c.f.~\cite[Proposition~3.2]{Mitrophanov2003}.
    % Furthermore if the chain $\MM$ obeys detailed balance (equivalently, is reversible), then $\|\mathbf{X}\|_2 \|\mathbf{X}^{-1}\|_2 = \sqrt{\max_iu_i/\min_iu_i} $, see~\cite[Theorem~3.2]{Mitrophanov2003}.

\end{remark}

% time dependent stuff -------------------
We conclude this section by 
showing how introducing a perturbation to $\MM$ can be understood as modifying the evolution of probabilities at (discrete) time $n=0$, see \cite{lucarini2025interpretableequationfreeresponsetheory}. In particular this argument will lead to an explicit expression (see \eqref{eq:assympt_linearresponse} below) for the linear response of the perturbed system in terms of the unperturbed one. Recall $\uu$ and $\vv$ are the unique invariant probability measure for $\MM$ and $\MM_\eps$ respectively, and set 
\begin{equation}\label{eq:wepsn}
    \ww_\eps(n) := \MM_\eps^n \uu = (\MM + \eps \PP)^n \uu, \quad n\geq 0. 
\end{equation}
Then, similarly to \eqref{eq:exp_mix},
\begin{equation*}
 \| \ww_\eps(n) - \vv\|_1 = \| \MM_\eps^n \uu - \vv \|_1 \leq \tau(\mathbf{M}_\eps)^n \|\uu - \vv\|_1
\end{equation*}
so, since $\tau(M_\eps)<1$ by Proposition~\ref{prop:ergodicity bound}, $\lim_{n\to \infty} \ww_\eps(n) =\vv$. 
Now, fixing $n\geq 0$ and expanding \eqref{eq:wepsn} in $\eps$, we get
\begin{equation}\label{eq:wepsn_expand}
    \ww_\eps(n) = \uu +\eps \ww_1(n) + R(\eps, n)
\end{equation}
where $\ww_1(n)$ is the linear response correction at time $n$, and  $R(\eps, n)\sim O(\eps^2)$ is the remainder, both null at time $n =0$. 
We now derive explicit formulae for $\ww_1(n)$ and $R(\eps, n)$.
By definition \eqref{eq:wepsn} and the expansion \eqref{eq:wepsn_expand}, for $n>0$ we can write
    \begin{align*}
        \ww_\eps(n) &= (\MM+ \eps \PP) \ww_\eps(n-1)
        = (\MM+ \eps \PP) (\uu +\eps \ww_1(n-1) + R(\eps, n-1))
        \\
        & = \uu + \eps(\MM \ww_1(n-1) + \PP\uu)  + \MM R(\eps,n-1) + \eps \PP R(\eps, n-1)+ \eps^2 \PP \ww_1(n-1),
    \end{align*}
    so that 
    \begin{align*}
        \ww_1(n) &= \MM \ww_1(n-1) + \PP \uu  \\ 
        R(\eps, n) &= (\MM + \eps \PP) R(\eps, n-1) + \eps^2 \PP \ww_1(n-1).
    \end{align*}
    Iterating these formulae, it is easy to see that 
    \begin{align}
        \ww_1(n) &= \sum_{k=0}^{n-1} \MM^k \PP \uu = \sum_{k=0}^{n-1} (\MM^k - \uu \mathbf{1}^{T}) \PP \uu \label{eq:linear_response_time_dependent}\\
        R(\eps, n)&= \eps^2 \sum_{k = 0}^{n-1} \MM_\eps^{n-1-k} \PP \ww_1(k) = \eps^2 \sum_{k = 0}^{n-1} (\MM_\eps^{n-1-k} - \vv \mathbf{1}^T) \PP \ww_1(k), \label{eq:R(eps,n)}
        % R(\eps, n+1) =& \MM_\eps^2 R(\eps, n-1) + \eps^2 \MM_\eps \PP \ww_1(n-1) + \eps^2 \PP \ww_1(n) \\
        %  \vdots&\\
        % =& \MM_\eps^{n+1}R(\eps, 0) + \eps^2 \sum_{k=0}^n \MM_\eps^{n-k}\PP \ww_1(k) \\
        %  =& \eps^2 \sum_{k = 0}^n \MM_\eps^{n-k} \PP \ww_1(k).
    \end{align}
    where we used assumption (C2) for $\PP$ and that $\ww_1(0) = R(\eps, 0) = 0$.
    Thanks to these expressions we can show that $\ww_1(n)$ converges exactly to the perturbed linear response $\left. \frac{d\vv}{d\eps}\right|_{\eps = 0} =: \vv_1$ for $n\to \infty$. 
    
    We start by noticing that, thanks to \eqref{eq:matrix_norm_bound},
    \begin{equation*}
        \|\ww_1(n)\| \leq \sum_{k = 0}^{n-1} \|\MM^k - \uu \mathbf{1}^{T}\| \|\PP\| \|\uu\| \leq 2\|\PP\| \sum_{k = 0}^{n-1} \tau(\MM)^k 
    \end{equation*}
    and, as $\tau(\MM)<1$, it follows 
    \begin{equation}\label{eq:unif_bound_w1}
        \|\ww_1(n)\| \leq \frac{2\|\PP\|}{1 - \tau(\MM)} \quad \text{for all } n\geq 0.
    \end{equation}
    This ensures that $\bar{\ww}_1 := \lim_{n\to \infty} \ww_1(n)$ is well defined.
    Using \eqref{eq:R(eps,n)},\eqref{eq:unif_bound_w1} and 
    \begin{equation*}
        \|\MM_\eps^k - \vv \mathbf{1}^T\|_1\leq 2 \tau(\MM_\eps)^k,
    \end{equation*}
    together with \cref{prop:ergodicity bound}, yields the uniform estimate 
     \begin{equation*}
       \|R(\eps, n)\| \leq \eps^2 \left( \frac{2\|\PP\|}{1 - \tau(\MM)} \right)^2 \quad \text{for all }n\geq 0. 
   \end{equation*}
    Taking $n\to \infty$ in \eqref{eq:wepsn_expand} therefore gives $\vv = \uu + \eps \bar{\ww}_1 + O(\eps^2)$ and hence 
    \begin{equation*}
       \bar{\ww}_1 =\left. \frac{d\vv}{d\eps}\right|_{\eps = 0} = \vv_1.
    \end{equation*}
    
    Next, to obtain an explicit expression for $\vv_1$, note that, since $\mathbf{1}^\top \PP = \mathbf{0}^\top$ by Assumption (C2), by \eqref{eq:linear_response_time_dependent} we can write
\begin{equation*}
    \vv_1 = \sum_{k = 0}^\infty (\MM^k-\uu\mathbf{1}^{\top}) \PP \uu = \left(\uu \mathbf{1}^\top + \sum_{k = 0}^\infty (\MM^k-\uu\mathbf{1}^{\top}) \right) \PP \uu
\end{equation*}
where the series is absolutely convergent since we have seen that
\begin{equation*}
    \sum_{k = 0}^\infty \|\MM^k-\uu\mathbf{1}^{\top} \| \leq \frac{2}{1 - \tau(\MM)}.
\end{equation*}
We define 
\begin{align}
   \GG := \uu\mathbf{1}^{\top} + \sum_{k=0}^\infty \left(\MM^k-\uu\mathbf{1}^{\top}\right)
\end{align}
so that $\vv_1 = \GG \PP \uu$. With some simple manipulations it is easy to see that 
\begin{equation}\label{eq:def:G}
   \GG = \sum_{k = 0}^\infty (\MM - \uu \mathbf{1}^\top)^k = \left(\mathbf{I}_N - \MM + \uu\mathbf{1}^{\top}\right)^{-1}
\end{equation}
where we have used the  definition of the Neumann series, and  where $\mathbf{I}_N$ is the identity $N\times N$ matrix. $\mathbf{G}$ is the generalized inverse of $\mathbf{I}_N - \MM$, which is not in general invertible, and it also called fundamental matrix of the Markov chain associated to $\MM$ \cite{seneta1993}.
Hence the stationary linear response $\vv_1$ can be written as
\begin{align}\label{eq:assympt_linearresponse}
    \vv_1 = \mathbf{G}\PP\uu.
\end{align}
 We note that the powers of $\MM$ do not converge to the null matrix, so the deflation through the projector matrix $\uu\mathbf{1}^{\top}$ is necessary for convergence. 
 
 By the same token, expectation values can be expanded in a linear response fashion. Let $\Psi$ denote a vector in $\R^N$, representing the coarse-graining of an observable function. The linear response correction to the expectation value of $\Psi$ with respect to $\ww_\eps(n)$, that we denote as $\langle \Psi \rangle_{\ww_\eps(n)}$, is:
\begin{equation}
    \frac{\mathrm{d}\langle \Psi \rangle_{\ww_\eps(n)}}{\mathrm{d}\eps}|_{\varepsilon = 0} = \Psi^{\top}\ww_1(n)
\end{equation}
By letting $n$ tend to infinity, i.e. replacing $\ww_1(n)$ by $\vv_1 = \GG \PP \uu$, we obtain the linear response at infinite time.

\section{Linear optimizations}\label{eq:optimizations}

The perturbative analysis of Markov chains presented in the previous section has served to analytically calculate sensitivity bounds of stationary measures \cite{seneta1993,Mitrophanov2003} but also to numerically compute linear responses of continuous-time dynamical systems \cite{Lucarini2016,SantosJSP}. In this section we address the question of finding perturbations of a given Markov matrix $\mathbf{M}$ that optimize (maximize/minimize) the value of functionals acting on matrices or vectors. Given a functional $f:\mathbb{R}^{N\times N}\longrightarrow \mathbb{R}$ acting on matrices, the maximization problem would read as:
\begin{equation}\label{eq:first_functional}
    \operatorname{arg max}_{\mathbf{P}\in \mathbb{R}^{N\times N}}\quad f(\mathbf{M}+\varepsilon \mathbf{P}),
\end{equation}
where $\varepsilon>0$ is a small parameter and where $\mathbf{P}$ is asked to satisfy Assumption~\ref{assum:perturb}.

The strategy carried out here is based on \cite{Antown2018}, where the authors exploit linear response formulae to transform optimizations over eigenproblems into linear ones. This way, the authors develop algorithms to optimize various physically meaningful quantities such as expectation values of observables and the rate of convergence to steady-state. The target of this section is to show how optimizing entropy-like functionals can be recast as linear-optimization problems which, in some cases, become analytically tractable.

\subsection{Linear response of entropy}\label{subsec:entropy}

The entropy of a strictly positive probability vector $\mathbf{u}$ is defined as:
\begin{equation}
    H(\mathbf{u}) = -\sum_{i=1}^Nu_i \log(u_i).
\end{equation}
A uniform probability vector represents a thermodynamic system where all microstates have the same probability of manifesting, so little information is known about the actual state of the system. This means that the system has high entropy, in fact, for the uniform distribution, maximum entropy. The opposite happens when the probability vector concentrates in a single microstate, yielding maximum information about the system, corresponding to  minimum entropy.

The entropy of the perturbed system, characterized by the perturbed invariant density $\mathbf{v}$--- solving Eq.~\eqref{eq:pert_eig}--- is to linear order in $\varepsilon$, $\mathbf{u}+\varepsilon\mathbf{v}_1$. Then, we can expand the entropy functional as:
\begin{subequations}
    \begin{align}
        H(\mathbf{u}+\varepsilon\mathbf{v}_1) &=- \sum_{i=1}^N(u_i + \varepsilon v_{1,i}) \log(u_i+\varepsilon v_{1,i}) = -\sum_{i=1}^N (u_i + \varepsilon v_{1,i})\left( \log(u_i) + \log\left(1+\varepsilon \frac{v_{1,i}}{u_i}\right) \right) = \\ &=H(\mathbf{u})-\varepsilon\sum_{i=1}^N \left( v_{1,i} \log(u_i)-v_{1,i}\right)+\mathcal{O}\left(\varepsilon^2\right)=H(\mathbf{u})-\varepsilon\sum_{i=1}^Nv_{1,i} \log(u_i)+\mathcal{O}\left(\varepsilon^2\right),
    \end{align}
\end{subequations}where we have used the fact that $ \log(1+\varepsilon \frac{v_{1,i}}{u_i}) = \varepsilon \frac{v_{1,i}}{u_i} + \mathcal{O}(\varepsilon^2)$ and that $\sum_i v_{1,i}=0$. The latter expression follows from the fact that both  $\mathbf{v}=\mathbf{u}+\varepsilon\mathbf{v}_1$ and $\mathbf{u}$ are probability vectors so $\sum_i v_i = \sum_i u_i = 1$.   
Thus, to optimize the entropy of the system, we search for a perturbation matrix $\mathbf{P}$ that maximizes the functional:
\begin{equation}\label{eq:max_entropy}
\frac{\mathrm{d}H(\mathbf{u}+\varepsilon\mathbf{v}_1)}{\mathrm{d}\varepsilon}|_{\varepsilon=0} =\mathbf{f}^{\top}\mathbf{G}\mathbf{P}\mathbf{u}  = \mathbf{f}^{\top} \mathbf{v}_1
\end{equation}
where $\mathrm{f}_i=\log(u_i)$ and $\mathbf{P}$ satisfies Assumption~\ref{assum:perturb}. This is consistent with the fact that the entropy $H(\mathbf{u})$ can be interpreted as the expectation value under the distribution $\mathbf{u}$ of the observable $\log(\mathbf{u})$, where the logarithm acts component-wise.

This optimization problem was solved in \cite[Section~4]{Antown2018} when maximizing the linear response of the expectation value of a generic observable. There the authors explicitly solve the Lagrange multiplier problem to find that the matrix $\mathbf{P}$ that maximizes Eq.~\eqref{eq:max_entropy} is given component-wise as:
\begin{equation}\label{eq:antown_expect}
    P_{ij} = \frac{u_j}{2\nu}\left( \sum_{k=1}^N\mathrm{f}_kG_{ki} - \frac{1}{m_j}\sum_{\ell:M_{\ell j}>0}\sum_{k=1}^N\mathrm{f}_kG_{k\ell} \right),
\end{equation}
where $m_j$ is the number of non-zero elements in the $j$th column of $\mathbf{M}$ and $\nu$ is a normalizing factor modulo a sign that depends on the Hessian matrix of the Lagrangian function. A similar derivation with another constraint is done later in the text in Section~\ref{eq:optimizations}.3, so more precise details are given then.

\begin{remark}
    Finite-time optimizations are possible if we substitute $\mathbf{G}$ by $\msout{\mathbf{G}(t)}$ \sout{defined below Eq.~\eqref{eq:linear_response_time_dependent}} by its finite truncation 
    \begin{equation*}
        \tilde{\GG}(n) := \uu\mathbf{1}^{\top} + \sum_{k=0}^{n-1} \left(\MM^k-\uu\mathbf{1}^{\top}\right) = \sum_{k=0}^{n-1} \left(\MM-\uu\mathbf{1}^{\top}\right)^k.
    \end{equation*}
     This way, we would obtain the perturbation $\mathbf{P}$ such that, for instance, entropy is maximized at time $n<\infty$. This finite-time approach is valid for the optimization of any functional including the ones presented here and in \cite{Antown2018}.
\end{remark}

\subsection{Kullback-Leibler divergence}\label{subsec:KL}

While optimizing entropy in absolute terms can be useful, it is also important to determine the perturbation that creates the most information relative to the unperturbed state as discussed in \cite{majda_abramov_book}. In other terms, we want to find the perturbation matrix $\mathbf{P}$ that maximizes the Kullback-Leibler (KL) divergence. The KL-divergence of the probability vector $\mathbf{v}$ with respect to $\uu$ is given by:
\begin{equation}\label{e9}
    \mathrm{D}_{\mathrm{KL}}(\vv,\mathbf{u}) = \sum_{i=1}^N v_i\log\left(\frac{v_i}{u_i}\right).
\end{equation}
Using the expression derived in Eq.\eqref{eq:assympt_linearresponse} for the response $\mathbf{v}_1$, we evaluate the KL-divergence between the leading order perturbed density $\mathbf{u} + \varepsilon \mathbf{v}_1$ and the unperturbed invariant vector $\mathbf{u}$ yielding--- c.f. \cite[Section~2.3]{majda_abramov_book}---:

\begin{subequations}\label{e10}
\begin{align}
    \mathrm{D}_{\mathrm{KL}}(\uu+\varepsilon \vv_1,\mathbf{u}) &= \sum_{i=1}^N (u_i+\varepsilon v_{1,i})\log\left( 1 +\varepsilon \frac{ v_{1,i}}{u_i} \right)= \\
    &= \sum_{i=1}^N (u_i+\varepsilon v_{1,i}) \left( \varepsilon \frac{v_{1,i}}{u_i} - \frac{\varepsilon^2}{2}\left( \frac{v_{1,i}}{u_i}\right)^2 + \mathcal{O}(\varepsilon^3)\right)= \\ 
    &=\frac{\varepsilon^2}{2}\sum_{i=1}^N\frac{v_{1,i}^2}{u_i} + \mathcal{O}\left(\varepsilon^3 \right) 
\end{align}
\end{subequations}

where we have used the second order Taylor expansion for $\log(1+x)=x - x^2/2 + \mathcal{O}(x^3)$ and the fact that $\sum_iv_{1,i}=0$, as already observed in the previous section. Maximizing the Kullback-Leibler divergence corresponds to find an optimal value for $\frac{1}{2}\sum_iv_{1,i}^2/u_i$. This term can be recast as a matrix product.   Indeed, define $\DD = (D_{ij})$ as a diagonal $N\times N $ matrix such that $D_{ii} = u_i>0$, for $i=1,\ldots,N$. Then, by expressing $\mathbf{v}_1=\mathbf{GPu}$ element-wise, we have
\begin{equation}\label{eq:rel_entr}
    \frac{1}{2}\sum_{i=1}^N \frac{1}{u_i} \left(\sum_{j=1}^N \sum_{\ell=1}^N  G_{i\ell}P_{\ell j}u_j\right)^2 = \frac{1}{2}\left\| \mathbf{D}^{-\frac{1}{2}}\mathbf{G}\mathbf{P}\mathbf{u} \right\|^2_2,
\end{equation}
where the matrix $\mathbf{D}^{-\frac{1}{2}}$ is a diagonal matrix with the $i$-th diagonal entry equal to $u_i^{-\frac{1}{2}}$ and such that $\mathbf{D}^{-\frac{1}{2}}\mathbf{D}^{-\frac{1}{2}} = \mathbf{D}^{-1}$.

The target is to find the matrix $\mathbf{P}$ that maximizes Eq.~\eqref{eq:rel_entr} subject to constraints (C1)-(C2)-(C3) in Assumption~\ref{assum:perturb}. The following similar problem was solved in \cite[Section 3.4]{Antown2018}:
\begin{equation}\label{eq:sect_3.4_antown}
    \operatorname{arg max}_{\mathbf{P}\in \mathbb{R}^{N\times N}} \left\| \mathbf{G}\mathbf{P}\mathbf{u} \right\|^2_2,
\end{equation}
subject to the same constraints as (C1)-(C2)-(C3) herein--- see also \cite[Eqs.~(34)-(37)]{Antown2018}---. The main difference between Eq.~\eqref{eq:sect_3.4_antown} and Eq.~\eqref{eq:rel_entr} is the presence of the matrix $\mathbf{D}^{-\frac{1}{2}}$ and the factor $\frac{1}{2}$. Hence, Algorithm~2 in \cite{Antown2018} is modified by replacing $\mathbf{G}$ in Eq.~\eqref{eq:sect_3.4_antown} by $2^{-\frac{1}{2}}\mathbf{D}^{-\frac{1}{2}}\mathbf{G}$. Furthermore, because of the quadratic structure of Eq.~\eqref{eq:rel_entr}--- and being of order $\mathcal{O}(\varepsilon^2)$---, if $\mathbf{P}$ maximizes relative entropy with respect to the unperturbed invariant measure, so will $-\mathbf{P}$.

\begin{remark}
   One can employ the same formulation as Eq.~\eqref{e10} to drive the system to having a stationary density closest to a prescribed distribution in the KL-divergence sense. Indeed, letting $\mathbf{w}=(w_i)$ denote a $N$-dimensional target probability vector, we want to minimize $\mathrm{D_{KL}}(\mathbf{u}+\varepsilon\mathbf{v}_1,\mathbf{w})$ using Eq.~\eqref{eq:rel_entr}, where, instead, the diagonal elements of $\DD$ are given by $D_{ii} = w_i$.
\end{remark}

\subsection{Entropy production}\label{subsec:Entropy_prod}
Entropy production is a paradigmatical observable to quantify how far physical systems are from detailed balance conditions due to the presence of non-conservative forces or the interaction with multiple heat baths at different temperatures. We refer the reader to \cite{schnakenberg} for a canonical reference on irreversible processes for Markov Chain systems. A non vanishing entropy production is related to irreversibility conditions and reflects the arrow of time as it can be defined in terms of the \textit{forward} and \textit{backward} trajectories in time \cite{Lebowitz_Spohn_1999}. Entropy production is a key quantity in stochastic thermodynamics where the classical macroscopic theory of thermodynamics is developed for mesoscopic and microscopic systems, see \cite{Seifert_2012}. Overall, entropy production  connects microscopic stochastic dynamics with macroscopic thermodynamic behavior, providing insight of the second law of thermodynamics at a probabilistic level \cite{Seifert_2012}.

The entropy production of a (column-normalized) Markov matrix $\MM$ is defined here as in \cite{schnakenberg}:
\begin{equation}\label{eq:entropy_production}
    s(\MM) \coloneqq \sum_{i=1}^N\sum_{j=1}^Nu_iM_{ji}\log\left( \frac{u_iM_{ji}}{u_jM_{ij}} \right),
\end{equation}
where, in general, $s(\MM)\geq 0$. The above definition makes sense if, when $M_{ij}>0$, then $M_{ji}>0$ for all $ i,j=1,\dots ,N$. When this assumption breaks, the system has genuine unidirectional transitions; it is an area of active research how to extend this definition to such systems \cite{Busiello_2020}. If there are no irreversible processes, the system obeys the detailed balance property $u_iM_{ji}/u_jM_{ij}=1$ for every $i,j=1,\ldots,N$ and the entropy production vanishes, $s(\MM)=0$. 
The question we want to answer is: what is the matrix $\PP$ such that the Markov matrix $\MM + \varepsilon\PP$ is closer to detailed-balance and has $\uu$ as an invariant probability vector? 

Consider the so-called additive \textit{reversibilization} of $\mathbf{M}$ (see e.g.~\cite{bremaud_markov_chains}),
\begin{equation*}
    \MM_r = \frac{1}{2}\left( \mathbf{M}+\mathbf{D}\mathbf{M}^{\top}\mathbf{D}^{-1} \right)
\end{equation*}
where $\DD$ is, as in the previous section, the diagonal matrix with $\uu$ on the diagonal. Since, by construction, $\MM_r$ satisfies detailed balance, a natural candidate for the perturbation bringing the matrix $\mathbf{M}$ closer to detailed balance would be  
\begin{equation}\label{eq:reversibilized_perturbation}
    \mathbf{P}_r =  \MM_r - \mathbf{M}.
\end{equation}
In this case, trivially, we then have that $(\mathbf{M}+\mathbf{P}_r)\mathbf{u} = \mathbf{u}$ and $s(\mathbf{M}+\mathbf{P}_r)=0$. The issue with this choice is that we aim at minimizing the distance to reversibility by \textit{small} perturbations to the original chain $\mathbf{M}$, while $\mathbf{P}_r$ might provoke significant changes in the chain as quantified by a large amplitude $\|\mathbf{P}_r\|_F\gg \varepsilon$. To solve such problem, we instead use the entropy production \eqref{eq:entropy_production} as proxy for reversibility and expand $s(\mathbf{M}+\varepsilon \mathbf{P})$ perturbatively in terms of $\varepsilon$ to obtain the leading order contribution of $\mathbf{P}$.

Because we aim at leaving the invariant measure unperturbed, $\uu$ must satisfy $\left(\MM + \varepsilon\PP \right)\uu = \uu$, and it follows that $\PP\uu = \mathbf{0}$. Hence, the matrix $\mathbf{P}$ in the minimization problem should satisfy another constraint on top of (C1)-(C2)-(C3):
\begin{align*}
\text{(C4)} &\quad \PP \uu = \mathbf{0} \quad \text{(preserve invariant probability)}
\end{align*}
To obtain the leading order contribution of $\mathbf{P}$ in powers of $\varepsilon$, we take the derivative of $s(\MM+\varepsilon\PP)$ with respect to $\varepsilon$ evaluated at $\varepsilon=0$. The main difficulty lies in expanding the logarithm inside Eq.~\eqref{eq:entropy_production} in powers of $\varepsilon$:
\begin{equation}
    \log \left(\frac{u_i(M_{ji} + \varepsilon P_{ji})}{u_j(M_{ij}+\varepsilon P_{ij})}\right) = \log\left( \frac{u_iM_{ji}}{u_jM_{ij}} \right) + \varepsilon\left( \frac{P_{ji}}{M_{ji}} - \frac{P_{ij}}{M_{ij}} \right) + \mathcal{O}\left(\varepsilon^2 \right).
\end{equation}
Thus, the leading order term of $s(\mathbf{M}+\varepsilon\mathbf{P})$ in powers of $\varepsilon$ reads as:

\begin{subequations}\label{eq:min_ent_prod}
\begin{align}
\frac{\mathrm{d}s(\MM+\varepsilon\PP)}{\mathrm{d}\varepsilon}|_{\varepsilon=0} &=    \sum_{i=1}^N\sum_{j=1}^Nu_i\left[P_{ji}\log\left( \frac{u_iM_{ji}}{u_jM_{ij}}\right) + P_{ji} - \frac{M_{ji}P_{ij}}{M_{ij}}\right] \\    
    &=\sum_{i=1}^N\sum_{j=1}^Nu_i\left[P_{ji}\log\left( \frac{u_iM_{ji}}{u_jM_{ij}}\right) - \frac{M_{ji}P_{ij}}{M_{ij}}\right] \\   
    &= \sum_{i=1}^N\sum_{j=1}^N\left[u_jP_{ij}\log\left( \frac{u_jM_{ij}}{u_iM_{ji}}\right)  -u_i \frac{M_{ji}P_{ij}}{M_{ij}}\right]  = \sum_{i=1}^N\sum_{j=1}^NC_{ij}P_{ij},
\end{align}
\end{subequations}
where we have used the fact that $\sum_{i=1}P_{ji}u_i =0 $ for every $j=1,\ldots,N$ and the matrix $\CC = (C_{ij})$ is a matrix whose entries are a function of those of $\MM$ and $\uu$:
\begin{equation}
     C_{ij} = u_j\log\left( \frac{u_jM_{ij}}{u_iM_{ji}}\right)  -u_i \frac{M_{ji}}{M_{ij}}.
\end{equation}
This way, the problem has the same linear and entry-wise structure as Eq.~\eqref{eq:max_entropy}. So we consider the problem of minimizing the linear functional of Eq.~\eqref{eq:min_ent_prod} subject to the constraints (C1)-(C2)-(C3)-(C4).
Because of the linearity and entry-wise structure of the problem, the sparsity constraint (C3) is automatically satisfied if the correct Lagrange multiplier equation is solved for the relevant $P_{ij}$, namely, where $M_{ij}>0$. To minimize Eq.~\eqref{eq:min_ent_prod}, we employ the Lagrange multipliers in the same way as \cite[Section~4]{Antown2018}, with the additional (C4) constraint. The corresponding Lagrangian function reads as:
\begin{equation}\label{eq:1}
    L(\PP,\mathbf{r},\nu)  = \sum_{i=1}^N\sum_{j=1}^NC_{ij}P_{ij} - \mathbf{r}^{\top}\PP^{\top} \mathbf{1} - \qq^{\top}\PP\uu - \nu \left(\|\PP \|_F^2 -1\right),
\end{equation}
where $\qq$, $\mathbf{r}$ and $\nu$ are the Lagrangian multipliers. To maximize this Lagrangian over the parameters, we need to compute the gradient with respect to each entry of $\PP$.  Because $\PP$ must have the same sparsity mask as $\MM$, one can impose the same sparsity mask to $\CC$ without changing the value of $\sum_{i=1}^N\sum_{j=1}^NC_{ij}P_{ij}$. We then compute the derivative of $L$ with respect to $P_{ij}$ and equalize it to zero:
\begin{equation}
    \frac{\partial L}{\partial P_{ij}}\left(\PP,\mathbf{r},\qq,\nu \right) = C_{ij}-r_j - q_iu_j-2\nu P_{ij} = 0,
\end{equation}
from where we obtain:
\begin{equation}\label{eq:4}
    r_j = C_{ij} -q_iu_j-2\nu P_{ij}. 
\end{equation}
Define $Z = \{ (i,j)\in \mathbb{N}^{2} : M_{ij}>0 \}$, and $Z_j^{c}=\{ i \in \mathbb{N}: (i,j) \in Z \}$. From the previous display, we find:
\begin{equation}\label{eq:2}
\sum_{i:(i,j)\in Z}r_j = |Z_j^c|r_j = \sum_{i:(i,j)\in Z}C_{ij} - u_j \sum_{i:(i,j)\in Z}q_i.
\end{equation}
This way we get a formula for $r_j$. We now plug the previous display into Eq.~\eqref{eq:4} to find the expression for $P_{ij}$:
\begin{equation}\label{eq:lagrange_matrix_formula}
    P_{ij} = \frac{1}{2\nu}\left[ C_{ij} - q_iu_j - \frac{1}{|Z_{j}^c|}\left( \sum_{i:(i,j)\in Z}C_{ij} - u_j \sum_{i:(i,j)\in Z}q_i \right) \right]
\end{equation}
we are left with finding $q_i$ and $\nu$.

We now multiply Eq.~\eqref{eq:1} by $u_j$, 
\begin{equation}
    C_{ij}u_j-r_ju_j - q_iu_j^2-2\nu P_{ij}u_j = 0
\end{equation}
and then we sum over $j=1,\ldots,N$, and apply the condition $\PP\uu = \mathbf{0}$, for each value of $k=1,\ldots,N$:
\begin{equation}
    \sum_{j:(k,j)\in Z}C_{kj}u_j - \sum_{j:(k,j)\in Z} r_ju_j - q_k\sum_{j:(k,j)\in Z}u_j^2=0.
\end{equation}
Next, we plug in the value of $r_j$ obtained in Eq.~\eqref{eq:2}:
\begin{equation}
    \sum_{j:(k,j)\in Z}C_{ij}u_j - \sum_{j:(k,j)\in Z}u_j \left[ \frac{1}{|Z_j^c|}\sum_{i:(i,j)\in Z}C_{ij} - \frac{u_j}{|Z_j^c|}\sum_{i:(i,j)\in Z}q_i \right] - q_k \sum_{j:(k,j)\in Z}u_j^2 = 0.
\end{equation}
Now we rearrange the previous equation into a linear equation on $q_{k}$ as:
\begin{equation}\label{eq:3}
\alpha_k + \sum_{j:(k,j)\in Z}^n\beta_j\sum_{i:(i,j)\in Z} q_i - \xi_k q_k =0
\end{equation}
where $\alpha_k$ and $\beta_j$ are defined as:
\begin{subequations}
    \begin{align}
        \alpha_k &=  \sum_{j:(k,j)\in Z}C_{ij}u_j - \sum_{j:(k,j)\in Z}  \frac{u_j}{|Z_j^c|}\sum_{i:(i,j)\in Z}C_{ij},  \\
        \beta_j &= \frac{u_j^2}{|Z_j^c|}, \\
        \xi_k &= \sum_{j:(k,j)\in Z}u_j^2.
    \end{align}
\end{subequations}

Equation~\eqref{eq:3} can be recast in the following form:
\begin{equation}\label{eq:linear_equation}
    \alpha_k + \sum_{i=1}^N \gamma_i q_i - \xi_k q_k =0,
\end{equation}
where $\sum_k \gamma_k = 1$. Notice that the associated matrix is invertible. Indeed, Eq.~\eqref{eq:linear_equation} can be written as:
\begin{equation}
    \left( \Xi + \mathbf{1}\boldsymbol{\gamma}^{\top} \right) \qq = \boldsymbol{\alpha}
\end{equation}
where $\Xi$ is a diagonal matrix with $\Xi_{ii}=\xi_i>0$, $\boldsymbol{\alpha} = (\alpha_i)$ and $\boldsymbol{\gamma}=(\gamma_i)$. Since $\Xi$ is invertible and $\mathbf{1}\boldsymbol{\gamma}^{\top}$ rank-one, the inverse of $\Xi + \mathbf{1}\boldsymbol{\gamma}^{\top}$ is given by the Sherman-Morrison formula:
\begin{equation}
    \left( \Xi + \mathbf{1}\boldsymbol{\gamma}^{\top} \right)^{-1} = \Xi^{-1} - \frac{\Xi^{-1}\mathbf{1}\boldsymbol{\gamma}^{\top}\Xi^{-1}}{1 + \boldsymbol{\gamma}\Xi^{-1}\mathbf{1}}.
\end{equation}

Finally $\nu$ is chosen so that the Frobenius norm constraint (C1) is satisfied and the sign is chosen negative so that the Hessian matrix of the Lagrangian function is negative; see \cite[Eq.~(57)]{Antown2018}. This guarantees that the matrix $\mathbf{P}$ in Eq.~\eqref{eq:lagrange_matrix_formula} is a minimizer.

An alternative, yet equivalent, approach can be taken by constructing a suitable vector space that contains the matrices satisfying constraints (C2) and (C4). For the particular problem of entropy production minimization, we run both methods on randomly initialized sparse and diagonal-dominant Markov matrices to illustrate such equivalence and also show that the methods, indeed, minimize the entropy production functional; see Appendix~\ref{sec:projection_space} for details.

% {\color{red}\paragraph{Special case: $M_{ij}>0$.} If $|Z| = n^2$, this means that $M_{ij}>0$ for every $i,j=1,\ldots,n$. Also, $|Z_j^c|=n$ for every $j=1,\ldots,n$. Hence, the resulting system of Eq.~\eqref{eq:3} reads as:
% \begin{equation}
%     \alpha_k +\frac{1}{n} \sum_{i=1}^n q_i - \| \uu \|_2^2 q_k =0,
% \end{equation}
% The resulting linear system for $q_i$ can be written in the same way as Eq.~\eqref{}, but in this case $\Xi = \|\mathbf{u}\|_2^2\mathbf{I}$ and $\gamma_i=1/n$ for every $i=1,\ldots,N$. This matrix is NOT invertible with $\ker(\AAA) = \mathrm{span}(\mathbf{1})$, but $\boldsymbol{\alpha}$ is in $\mathrm{rank}(\AAA)$, so a solution exists of the form:
% \begin{equation}
%     \mathbf{q} = \AAA^{\dagger}\boldsymbol{\alpha} + \xi \mathbf{1},
% \end{equation}
% for any $\xi$ in $\mathbb{R}$. Therefore there is no unique solution. However, we can take $\xi$ so that $\sum_iq_i = n\cdot \xi,$ this means that $\qq \perp \ker(\AAA)$.

% Actually, in this case, a solution to Eq.~\eqref{eq:6} is also given by a solution to the following problem:
% \begin{equation}
%     \left(A - \frac{1}{n}\mathbf{1}\mathbf{1}^{\top} \right)\qq = I\qq = \boldsymbol{\alpha},
% \end{equation}
% hence, $\qq = \boldsymbol{\alpha}$.
% }

\begin{remark}
In the Section~5 of \cite{Antown2018}, the authors demonstrate that the linear correction to the absolute value of largest nonunit eigenvalue of $\mathbf{M}$ can be written in the same entry-wise structure as Eq.~\eqref{eq:min_ent_prod}, see \cite[Eq.~(68)]{Antown2018}. This means that the formula Eq.~\eqref{eq:lagrange_matrix_formula} is readily applicable in such problem, with the extra constraint of invariant vector preservation. This is of interest for problems where one aims at accelerating/decelerating mixing or convergence to stationarity while preserving the same invariant measure \cite{pavliotisbook2014}.
\end{remark}

\section{From flows to Markov matrices and back: the Ulam approach}\label{sec:ulam}

As mentioned in the introduction, Markov matrices arise naturally in the study of continuous and discrete time dynamical systems and stochastic differential equations (SDEs). In this section we briefly review a standard methodology of obtaining Markov matrix representations of flows from time series and we provide a strategy to numerically reconstruct a vector field from perturbations of Markov matrices.

Let us consider a $d$-dimensional It\^o SDE with drift $\FF : \R^d \longrightarrow \R^d$, and volatility matrix $\Sigma \in \R^{d\times d}$:
 \begin{equation}\label{eq:sto ode 2}
\dd \xx  = \FF (\xx)\dd t+\Sigma(\xx)\dd \mathbf{W}_{t},
\end{equation}
where $\mathbf{W}_t$ is a $d$-dimensional Wiener process. An equivalent representation of Eq.~\eqref{eq:sto ode 2} is provided by the evolution law of probability density functions in $\mathbb{R}^d$. This is captured by the Fokker-Planck equation \cite{risken}:
\begin{equation}\label{eq:fpe 2}
\partial _t\rho = \LLL \rho= -\nabla \cdot \lp \FF\rho \rp   +\frac{1}{2}\nabla^2:\lp  \Sigma \Sigma^{\top}\rho\rp ,
\end{equation}
where the symbol ``$:$'' is the Hadamard, entry-wise product. 
The unique strong solution $(X_t)_{t\geq 0}$ of~\eqref{eq:sto ode 2} defines a
Markov semigroup $\{P_t\}_{t\geq 0}$ on $L^2_{\rho_0}$ where $\rho_0$ is the invariant density of the system, via
$(P_t f)(x) = \mathbb{E}^x[f(X_t)]$; the family $\{e^{t\mathcal{L}}\}_{t\geq 0}$
is its $L^2_{\rho_0}$-adjoint, and is a $C_0$-semigroup by invariance of
$\rho_0$~\cite{engel2000}.

Having access to the operator $\LLL$ provides essential information on the system in Eq.~\eqref{eq:sto ode 2}. For instance, the rate at which correlations decay, its response and statistics are determined by such an operator \cite{baladi,chekroun2019c,Santosgutierrez_2022}. However, tracing the Fokker-Planck operator $\LLL$ analytically is generally undoable, apart from the most elementary cases: see e.g. \cite{metafunes2002}. In general, the spectral properties of $\LLL$ are numerically estimated from the analysis of time series obtained from the system in question, regardless of it being stochastic or deterministic. Indeed, the so-called Ulam's method grants a way to approximate the semigroup $\{ e^{t\LLL}\}_{t\geq0}$ via Markov chains \cite{ulam1964,froylandapproximating1998}. To do this, we subdivide the $d$-dimensional phase space $\mathcal{X}\subseteq \R^{d}$ into $N$ non-intersecting \emph{cells} or \emph{boxes} $\{B_i\}_{i=1}^{N}$ and we define $\chi_{B_i}$ as the indicator function on the box $B_i \subset \mathcal{X}$. We denote the $N$-dimensional vector space spanned by this set of functions as $\mathcal{V}_N$.  We approximate the continuous time interval as a sequence of discrete times $\lbrace j\tau, \, j\in \mathbb{N}\rbrace$ with time step $\tau>0$, also called transition time. Denoting $\mathcal{P}_N$ as the conditional expectation onto $\mathcal{V}_N$, the projected exponential operator $\mathcal{P}_Ne^{\tau\LLL }: \mathcal{V}_N \longrightarrow \mathcal{V}_N$, hence, admits a matrix representation $\mathbf{M}_{\tau}$, where each element is given by:
\begin{equation}\label{projected transfer operator}
M_{\tau,ij}:=\frac{1}{\rho_0(B_i)}\int_{B_i}e^{\tau \LLL}\chi_{B_j}(\xx)\rho_0(\dd \xx),
\end{equation}
with $i,j=1,\ldots,N$. While Ulam's method has proven to be effective, it remains to fully understand the nature of these approximations as the resolution $N$ increases in deterministic systems. In this case, rigorous results are restricted to one-dimensional systems \cite{li1976} and those systems with smooth invariant measures \cite{froylandapproximating1998}. For stochastic systems, uniform ellipticity in the diffusion law and a confining drift term guarantee that the operator $e^{\tau\mathcal{L}}$ is compact in $L^2_{\rho_0}$. In this case, the eigenpairs of Ulam matrices $\mathbf{M}_{\tau}$ converge to the true values if the projection $\mathcal{P}_N$ converges to the identity in the operator norm as $N$ tends to infinity--- see \cite[Theorems 5 and 6]{osborn1975} for convergence results for eigenvectors and eigenvalues---.

At this stage, formula Eq.~\eqref{projected transfer operator} links SDEs with the Markov matrix framework in Section~\ref{sec:pert_markov}. Once an adequate Markov matrix is constructed--- for a fixed $\tau$--- as an approximation to the Fokker-Planck semigroup, all the optimization algorithms listed in Section~\ref{eq:optimizations} are readily applicable. In general, the validity of Ulam's approximation depends on the transition time $\tau$ used and the implementation of the matrix logarithm \cite{KLUS2020_GEDMD}. Broadly, the coarser the gridding the longer $\tau$ should be to avoid the problem of excessive artificial diffusion \cite{generatorfroyland,Tantet2018}. We refer the reader to \cite{METZNER2007} for further details on optimized methods for generator estimation with Ulam-like basis functions. From now on, the subscript $\tau$ that indicates the transition time in Eq.~\eqref{projected transfer operator} will be omitted when it is not relevant to the discussion.

\subsection{Phase space reconstruction of the perturbation operator}\label{subsec:drift_reconstr}

Once the optimal perturbation matrix $\mathbf{P}$ is found, it is important to reconstruct the corresponding perturbation at the SDE level. In other words, we want to find the vector field such that adding it to Eq.~\eqref{eq:sto ode 2} corresponds to the perturbed Markov chain system. In this section, we closely follow \cite{KLUS2020_GEDMD} to perform \textit{system identification}, that is, estimating the drift and the diffusion terms of stochastic differential equations from data using transfer operator techniques. First, we briefly summarize their approach. Then, we will describe how this method can be used to estimate the vector field associated with the optimal perturbation $\mathbf{P}$. 

We introduce the Koopman operator (also known as the Markov semigroup in stochastic settings) acting on observables $f$ as
\begin{equation}
    (\mathcal{K}^tf)(\mathbf{x}) = \mathbb{E}[f(\mathbf{x}_t)|\mathbf{x}_0=\mathbf{x}],
\end{equation}
where $\mathbf{x}_t$ is the stochastic process generated by \eqref{eq:sto ode 2} conditioned on starting at point $\mathbf{x}$, and $\mathbb{E}[\cdot]$ is the expected value over all noise realizations. The generator of the Koopman operator $\mathcal{L}^*$ is defined as $\mathcal{K}^t=e^{\mathcal{L}^*t}$ and is the adjoint of the Fokker Planck operator $\mathcal{L}$
\begin{equation}
    \mathcal{L}^* = \mathbf{F}\cdot \nabla  + \frac{1}{2}\Sigma\Sigma^\top : \nabla^2.
\end{equation}

There are multiple algorithms to evaluate the spectral properties of the Koopman operator and its generator \cite{COLBROOK_Multiverse,KLUS2020_GEDMD}. Of particular interest here is the fact that the standard Ulam method represents a specific case of such algorithms and can be used to estimate the properties of the Koopman operator and its generator \cite{Klus_2016,klus_2018_approx_transfer_operator}. In particular, a generalized eigenvalue problem for the matrix $ \log (\mathbf{M})/  \tau$ provides an estimate of the eigenfunctions $\varphi_i$ and eigenvalues $\lambda_i$ of the Koopman generator $\mathcal{L}^* \varphi_i = \lambda_i \varphi_i$. 

\textit{System identification} can then be performed by considering the action of Koopman generator on the state observable $f(\mathbf{x}) = x_i$. The idea is that the $i$-th element of the vector field $\mathbf{F}$ is related to the Koopman generator and its spectrum in the following manner: 
\begin{equation}\label{eq:spectral_reconstruction}
    F_i(\mathbf{x})=\mathcal{L}^* x_i = \mathcal{L}^*\left(\sum_{k=1}^{K}  \xi_{ik} \varphi_k(\mathbf{x})\right) = \sum_{k=1}^K \lambda_k \xi_{ik} \varphi_k(\mathbf{x}),
\end{equation}
where $K$ indicates the number of resolved eigenfunctions and the matrix $\Xi = (\xi_{ik})$ contains the Koopman modes, that is, the decomposition of the state variable $x_i$ in terms of the $K$ Koopman eigenfunctions. Repeating the operation for all the variables of the system gives the full vector field.

The diffusion matrix $\Sigma\Sigma^\top$ associated with \eqref{eq:sto ode 2} can then be approximated in a similar way by considering the action of the generator on the quadratic observable $f(\mathbf{x})=x_ix_j$ \cite{KLUS2020_GEDMD}:
\begin{equation}
    \mathcal{L}^*( x_ix_j) = F_i(\mathbf{x})x_j + F_j(\mathbf{x})x_i + (\Sigma\Sigma^\top)_{ij}.
\end{equation}
This methodology allows one to estimate the drift and diffusion matrix of the stochastic dynamics given its Markov chain representation through the Ulam method. Indeed, once the spectrum of $\mathbf{M}$ is computed, we can seek for a decomposition of generic observables in the eigenbasis. Given an observable $f$, we can first estimate its decomposition into the Ulam characteristic functions as $f(\mathbf{x}) = \sum_j \hat{f}_j \chi_{B_j}(\mathbf{x})$. Once the coefficients $\hat{f}_j$ are known, the same generalized eigenvalue problem mentioned above allows us to write the decomposition of the observable in to the Koopman eigenfunctions $f(\mathbf{x})=\sum_k f_k \varphi_k(\mathbf{x})$, see \cite{Klus_2016} for details. This way, choosing the observable $f(\mathbf{x})=x_i$, one can readily apply Eq.~\eqref{eq:spectral_reconstruction}.

We propose to use this approach to estimate the drift and diffusion associated to the optimal perturbation $\mathbf{P}$. We denote as $\widehat{\mathbf{F}}$ and $\widehat{\mathbf{F}}_\varepsilon$ the drift estimate of the Markov matrix $\mathbf{M}$ and $\mathbf{M} +\varepsilon\mathbf{P}$ respectively, where $\varepsilon$ is chosen to be in a linear response regime. The drift associated with the perturbation $\mathbf{P}$ is estimated as $(\widehat{\mathbf{F}}_\varepsilon - \widehat{\mathbf{F}}) / \varepsilon$. The diffusion matrix associated with the perturbation can be estimated in a similar way. 

\begin{remark}
    The relation $F_i(\mathbf{x})=\mathcal{L}^* x_i$, for every $i=1,\ldots,d$, can also be used in a finite difference perspective, without resorting to computing the spectrum of $\mathbf{M}$ or $\mathbf{M}+\varepsilon\mathbf{P}$. The entries of a perturbation matrix $\mathbf{P}$ represent the infinitesimal modification of probability rate fluxes out of each box in the discretization. Hence, the mean velocity is the sum over all these possible transitions weighted by their rates, given by $\mathbf{P}$, and the displacement vectors. With this ansatz, we get a finite-difference approximation formula for the vector field:
\begin{equation}\label{eq:drift_recons}
\left[\left(\widehat{\mathbf{F}}_\varepsilon - \widehat{\mathbf{F}}\right)(\mathbf{c}_i) \right]_k \approx \sum_{j=1}^N P_{ji} \left(c^{(k)}_j-c^{(k)}_i\right),    
\end{equation}
where $\left[\left(\widehat{\mathbf{F}}_\varepsilon - \widehat{\mathbf{F}}\right)(\mathbf{c}_i) \right]_k$ is the $k$-th component of the vector field $\widehat{\mathbf{F}}_\varepsilon - \widehat{\mathbf{F}}$ evaluated at box-center $\mathbf{c}_i$. In Appendix~\ref{app:finite_difference} provide a perturbative derivation of Eq.~\eqref{eq:drift_recons}. 
\end{remark}

It remains to explain how the estimated perturbation vector field $(\widehat{\mathbf{F}}_\varepsilon - \widehat{\mathbf{F}}) / \varepsilon$ depends on the transition time $\tau.$ For different values of $\tau$, the resulting transition matrices $\mathbf{M}$ are different and, hence, the perturbation matrix $\mathbf{P}$ can also change, apart from its Frobenius norm which is fixed to unity by Assumption~\ref{assum:perturb}. The reason is that $\mathbf{P}$ is constructed in terms of the nonzero entries of $\mathbf{M}$, which may alter as a function of $\tau$. Indeed, as $\tau$ increases the matrix $\mathbf{M}$ becomes less and less sparse. Conversely, as $\tau$ approaches zero, transitions between boxes become rarer and $\mathbf{M}$ eventually converges to the identity matrix. We can then consider $\mathbf{P}$ is independent of $\tau$ if its evaluation on the unperturbed invariant measure, namely, $\mathbf{P}\mathbf{u}$, does not change with $\tau$. More importantly, the linear response defined in Eq.~\eqref{eq:assympt_linearresponse} should stay invariant with $\tau$. Note, however, that such a linear response formula is the limit of a time average for discrete times with unit time increment--- see Eq.~\eqref{eq:linear_response_time_dependent}---. In the Ulam perspective, the time discretization unit is $\tau$ and, hence, time averages are weighted by $\tau$.
% Analytical solutions to this problem are expected to be untractable, so here a numerical algorithm is described. First of all, we want to investigate whether we can reconstruct the vector field $\FF$ in Eq.~\eqref{eq:sto ode 2} from the Markov matrix discretization, $\mathbf{M}$, of $e^{\tau\LLL}$.

% The generator $\LLL$ of the semigroup is approximated by the logarithm of the Markov matrix normalized by the transition time $\tau$--- as follows from the Spectral Mapping Theorem \cite{engel2000}---:
% \begin{equation}\label{eq:matrix_log}
%     \LLL \approx \mathbf{L} = \frac{1}{\tau} \log\left(\MM \right) = \sum_{k=1}^{\infty} \frac{(-1)^k}{k}\left(\MM - \mathbf{I}_N \right)^k.
% \end{equation}
% Because $\MM$ is column-stochastic, the columns of $\mathbf{L}$ must add up to zero. The entries of the matrix $\mathbf{L}$ represent the rate of probability fluxes out of each box in the discretization. Hence, the mean velocity is the sum over all these possible transitions weighted by their rates, given by $\mathbf{L}$, and displacement vectors. With this ansatz, we get the approximation formula for the vector field:
% \begin{equation}\label{eq:drift_recons}
% \left[\FF(\mathbf{c}_i) \right]_k \approx \sum_{j=1}^N L_{ji} \left(c^{(k)}_j-c^{(k)}_i\right),    
% \end{equation}
% where $\left[\FF(\mathbf{c}_i) \right]_k$ is the $k$-th component of the vectorfield $\mathbf{F}$ evaluated at box-center $\mathbf{c}_i$.

\subsubsection{Drift-diffusion reconstruction in a two-dimensional double well}

To illustrate the validity of this method, we consider a paradigmatic set of equations representing equilibrium dynamics of a particle in a potential with two different competing attractors embedded in thermal (Gaussian) noise, namely, a double well potential.

The double well potential in two dimensions is a simple system displaying non-Gaussian statistics and multimodality. We consider the following SDE on $\mathbb{R}^2$:
\begin{align}\label{eq:dwp_symm}
    \dd \mathbf{x} = \mathbf{F}(x,y)+ \sigma \dd \mathbf{W}_t=\binom{x - x^3}{-y}\dd t + \sigma \dd \mathbf{W}_t,
\end{align}
where $\sigma>0$ and $\mathbf{W}_t=(W_t^{(1)},W_t^{(2)})^{\top}$, $i=1,2$, denotes two dimensional independent standard Wiener process. It is easy to see that the equation is driven by the potential $-V(x,y) = x^2/2 - x^4/4 - y^2/2$. Here, we select a noise strength of $\sigma = 0.4$, so that the mean time to switch between metastable states is in the order of 10 time units.

In order to approximate the Fokker-Planck semigroup, we integrated Eq.~\eqref{eq:dwp_symm} for $10^7$ time units with a timestep of $10^{-2}$ after removing a transient of 10 time units which started from the origin. Snapshots are recorded every $0.1$ time units. This provides a time series $\{(x_k,y_k) \}_{k=1}^{K}$, where $K = 10^8$. The selected spatial domain is $D =[-2,2]\times [-1.5,1.5]$ and it is subdivided into $2^{14}$ boxes $B_i$, where each axis is divided into $2^7$ segments. We estimate the Markov matrix from Eq.~\eqref{projected transfer operator}:
\begin{equation}\label{eq:transition_matrix_1D}
    M_{ij} = \frac{\#\Big\{ \{ (x_k,y_k) \}_{k=1}^{K-\tau/\mathrm{d}t} \in B_j  \land \{(x_k,y_k)\}_{k=1+\tau/\mathrm{d}t}^K\in B_i \Big\}}{\# \Big\{ \{(x_k,y_k)\}_{k=1}^{K-\tau/\mathrm{d}t} \in B_j\Big\}},
\end{equation}
where $\tau$ is the transition time. For the searching algorithm we employ the binary tree structure from GAIO \cite{gaio}. Two separate matrices have been constructed with two transition time values: $\tau = 0.1$ and $\tau = 0.5$. For these values, correlation functions of the system can be well approximated and the mean exit times between potential wells are obtained from the first non-unit eigenvalue of $\mathbf{M}$; see \cite{banerjee2026datadrivenanalysismetastabilitystochastic}. This suggests that the spectrum of $\mathbf{M}$ contains meaningful information about the system and we can employ it to perform systems identification following Section~\ref{subsec:drift_reconstr}. Indeed, the drift and diffusion of Eq.~\eqref{eq:dwp_symm} are well reconstructed and the results are shown in Appendix~\ref{app:resolution}, where we also plotted the results for a coarser resolution involving $2^{12}$ boxes. Boxes that are not populated by the stochastic trajectory are removed. Here, we only show the results for $\tau=0.1$, since the results of $\tau=0.5$ are the same.

\begin{figure}[H]
\centering
\includegraphics[scale=0.29]{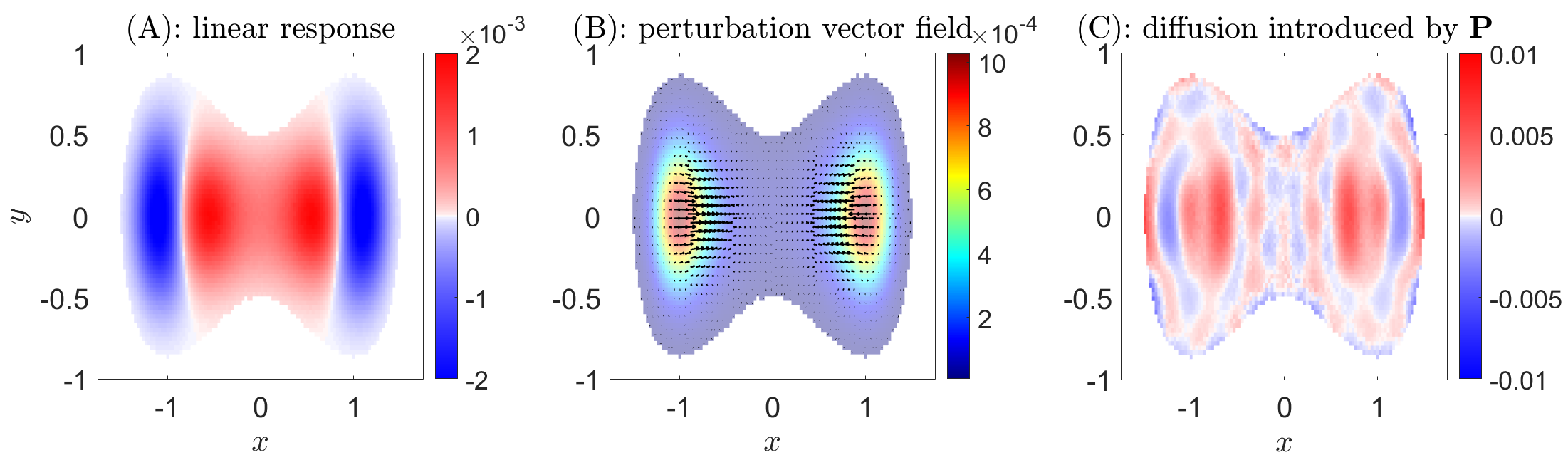}
\caption{\label{fig:dwp1_depth14}\textbf{Expectation value maximization: perturbation and linear response.} Panel (A) shows the linear response with respect to the perturbation $\mathbf{P}$. Panel (B) shows the perturbation vector field associated with the matrix perturbation $\mathbf{P}$ that maximizes the expectation value of $\Psi$. Such a vector field is superimposed on the invariant measure of the unperturbed system.  Panel (C) is the difference between the diagonal diffusion tensors obtained from $\mathbf{M}+\varepsilon \mathbf{P}$ and $\mathbf{M}$.}
\end{figure}

The invariant probability density of the system is approximated by the eigenvector $\uu$ that solves Eq.~\eqref{eq:inv_mes}, namely $\MM \uu = \uu$. Such a vector $\uu$ is shown underlying the arrow-plot in Fig.~\ref{fig:dwp1_depth14}(B). As mentioned earlier, since the system is in equilibrium, to all practical purposes the estimated Markov matrix is reversible, implying that the largest (in absolute value) one hundred and forty five eigenvalues were strictly real. These were used for the spectral reconstruction of the drift and diffusion.

% We tackle the problem of reconstructing the vector field $\mathbf{F}$ in Eq.~\eqref{eq:dwp_symm}, using the spectral procedure explained in Section~\ref{subsec:drift_reconstr}. The eigenfunction of the Fokker-Planck semigroup generator, necessary to implement the algorithm, are obtained from the eigenvector decomposition of the matrix $\mathbf{M}$ in Eq.~\eqref{eq:transition_matrix_1D}. Panels (A) and (B) in Fig.~\ref{fig:dwp1_depth14} show the absolute error between the analytical form of $\mathbf{F}$ and its reconstructed counterpart for the $x$ and $y$-directions respectively. We observe that while the areas where the invariant measure is larger, the errors are minimal, whereas areas at the border of the domain contain most of the errors. We also retrieved the diffusion tensor entailed by the matrix $\mathbf{M}$ and this is captured in panels (C) and (D) in Fig.~\ref{fig:dwp1}, where we only show the first diagonal component $\left(\Sigma\Sigma^{\top}\right)_{11}(x,y)$ and the off-diagonal component $\left(\Sigma\Sigma^{\top}\right)_{12}(x,y)$. From the model formulation, we know that $\left(\Sigma\Sigma^{\top}\right)_{11}$ and $\left(\Sigma\Sigma^{\top}\right)_{12}$ should be constant equal to $0.16$ and $0$, which is what the results indicate.

For illustration, we consider the maximization of the expectation value of an observable as in \cite{antown2022optimal} via \eqref{eq:antown_expect}. We consider the observable 
\begin{equation}\label{eq:observable}
    \Psi(x,y) = \left( 1 + (x/0.5)^2 + (y/0.4)^2\right)^{-4}. 
\end{equation}
This observable is a plateau-type of function that concentrates around the origin in the phase-space precisely where the invariant measure has a local minimum half-way between the two wells. The idea for choosing this observable is to examine to which extent the leading order response will push probability towards the origin in detriment of the two wells containing most of the probability. We now let $\mathbf{P}$ denote the optimal perturbation for \eqref{eq:antown_expect}. One can guess that in order to maximize the expectation value of $\Psi$, the perturbation $\mathbf{P}$ should transfer probability toward the origin. A value of $\varepsilon = 0.01$ is chosen. Although the results remain the same for other values of $\varepsilon$--- see Appendix~\ref{app:dependence_epsilon}---.

In Fig.~\ref{fig:dwp1_depth14} we summarize the results. Panel (A) shows the linear response of the system, this is, the leading order change in probability that maximizes the expectation value of $\Psi$. Panel (B) shows the reconstructed perturbation vector field (shown with arrows) superimposed to the unperturbed invariant measure $\uu$.  Lastly, in panel (C), we show that the difference between the first diagonal component of the diffusion tensor $\left(\Sigma\Sigma^{\top}\right)_{11}$ computed from $\mathbf{M}+\varepsilon \mathbf{P}$ and $\mathbf{M}$. We observe that the reconstructed perturbation vector field does precisely what one would expect, this is, to transfer probability toward the center of the domain where the observable $\Psi$ is mostly supported. The resulting net effects are captured by the linear response. A close observation of Panel (C), reveals that the magnitude of the diffusion tensor can reach values close to $7\%$ of the original diffusion value of $\sigma^2$, particularly around the wells. That, together with the patterns revealed by the plot--- preserved also in coarser resolution; see Appendix~\ref{app:resolution}--- suggest that the diffusion contained in $\mathbf{P}$ can be significant.

\section{Numerical experiments}\label{sec:numerics}
In this section, we assess the proposed framework through a range of numerical examples chosen to explore different aspects of the methodology.
Rather than applying the same optimization problem to every system, each example highlights particular features of the approach, so taken together, the experiments show how the framework can be adapted to different objectives, dynamical settings, and levels of model complexity.

For the Lanford map, we consider the maximization of the entropy functional and of the KL divergence, as described in Sections~\ref{subsec:entropy}--\ref{subsec:KL}. We then study a double-well potential, but his time with a rotational component, where we aim at minimizing entropy production. For the more complex Lorenz--63 system, we consider entropy functional maximization, as with the Lanford-map example. Then, we repeat the Lorenz--63 experiment but projected in 2D, to highlight the role of the loss of Markovianity. Finally, we apply the framework to a reduced Markov chain description of a turbulent flow based on unstable periodic orbits, identifying the perturbation that maximizes the average kinetic energy. We refer to Section~\ref{subsec:upo} for full details of the construction.

\subsection{The Lanford map}

In this section we analyze a stochastically perturbed Lanford map, which has also been investigated in the work of \cite{Antown2018} using the Markov matrix framework. Here, instead, we employ the algorithms presented in Section~\ref{eq:optimizations} to analyze its response to perturbations that maximize, 1) the linear response of the entropy functional and, 2) the KL-divergence of the perturbed measure with respect to the unperturbed state. For each case, the linearity of the response and its convergence to its asymptotic value is investigated.

\begin{figure}[H]
\centering
\includegraphics[scale=0.29]{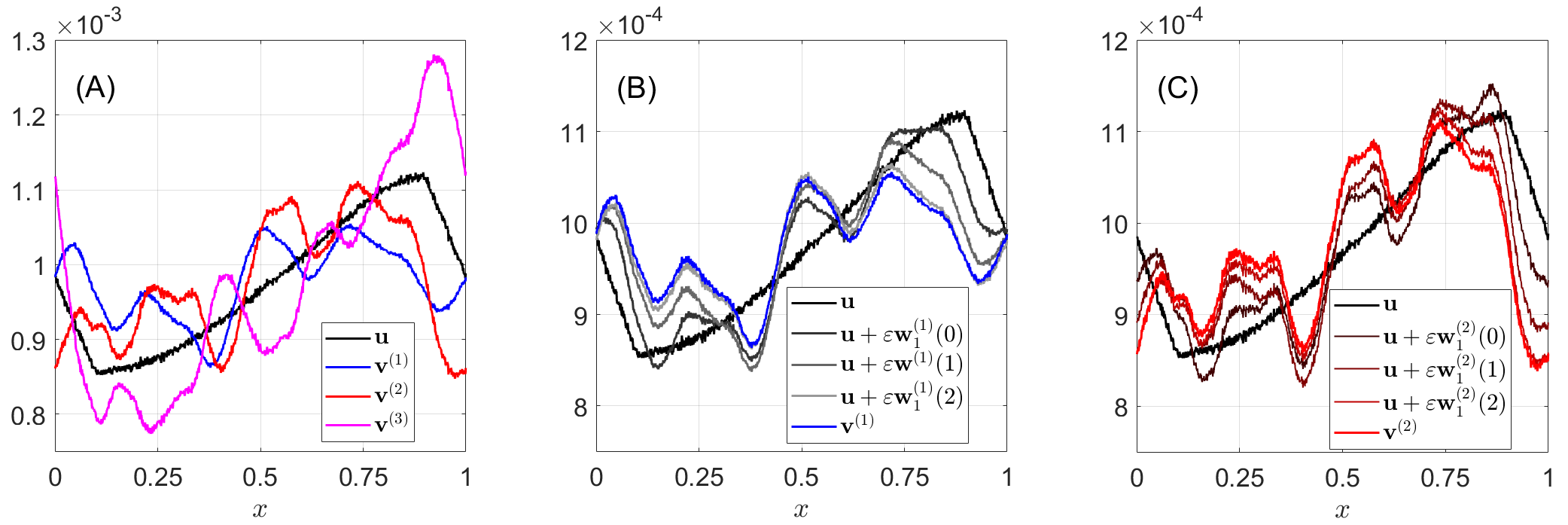}
\caption{\label{fig:lanford}\textbf{Perturbations of the Lanford map.} Panel (A): in black we show the invariant vector of the matrix in Eq.~\eqref{eq:transition_matrix}, in blue, red and magenta, the perturbed invariant vectors--- see Eq.~\eqref{eq:eig_prob}--- of the perturbation problems considered. Panel (B): Case of maximizing the entropy functional. In black we show the unperturbed invariant vector, in gray scales, the successive transient linear responses--- for the times indicated in the legend--- and in blue the perturbed invariant vector. Panel (B): same as Panel (A), but for the maximization of KL-divergence between perturbed and unperturbed invariant vectors.}
\end{figure}

The stochastic Lanford map $\mathfrak{L}:[0,1]\longrightarrow[0,1]$ is defined as:
\begin{equation}\label{lanford}
\mathfrak{L} (x) = 2x + \frac{1}{2}x(1-x) + \xi \mod 1    
\end{equation}
where $\xi$ is uniformly distributed on the interval $[-0.1,0.1]$. The Markov discretization is done as follows. First, the interval $[0,1]$ is discretized into $1024$ boxes $\{ B_i \}_{i=1}^{1024}$. 
Secondly, Eq.~\eqref{lanford} is integrated for $K= 10^8$ time steps producing a timeseries $\{x_k\}_{k=1}^{K}$, where $x_k$ is in $[0,1]$ and $x_1 = 0$. Third, the Markov matrix $\MM = (M_{ij})$ describing the pushforward map--- the discrete-time version of Eq.~\eqref{projected transfer operator}--- is constructed by:
\begin{equation}\label{eq:transition_matrix}
    M_{ij} = \frac{\#\Big\{ \{ x_k \}_{k=1}^{K-1} \in B_j  \land \{x_k\}_{k=2}^K\in B_i \Big\}}{\# \Big\{ \{x_k\}_{k=1}^{K-1} \in B_j\Big\}}.
\end{equation}
The resulting matrix $\MM$ has a unique stationary probability vector $\uu$ and it is plotted in the black line of Fig.~\ref{fig:lanford}(A).

The linear response of the entropy function is maximized using formula~\eqref{eq:antown_expect} to obtain a perturbation matrix $\mathbf{P}_1$. To maximize the KL-divergence from $\mathbf{u}$, the functional Eq.~\eqref{eq:rel_entr} is maximized by a matrix $\mathbf{P}_2$. We recall that the matrix $-\mathbf{P}_2$ would also maximize the KL-divergence from $\mathbf{u}$ due to the quadratic nature of the functional; see Sect.~\ref{subsec:KL}. For completeness we also analyze $\mathbf{P}_3=-\mathbf{P}_2$. Then, the perturbation size of $\varepsilon = 0.1$ is chosen so that $\MM + \varepsilon \PP_i$ is a (perturbed) Markov matrix with perturbed invariant measure $\mathbf{v}^{(i)}$ satisfying:
\begin{equation}\label{eq:eig_prob}
    \left( \MM + \varepsilon \PP_i \right)\mathbf{v}^{(i)}=\mathbf{v}^{(i)},
\end{equation}
for $i=1,2,3$. The respective perturbed steady states are plotted in blue and red in Fig.~\ref{fig:lanford}(A). At time $n=0$ the perturbation $\varepsilon\mathbf{P}_i$--- for $i=1,2,3$--- is applied to $\mathbf{u}$ and the perturbed system evolves in time towards its new steady state 
\begin{equation}
    \left( \MM + \varepsilon \PP_i \right)^{n}\mathbf{u} = \uu + \varepsilon \mathbf{w}_1^{(i)}(n) + \mathcal{O}\left(\varepsilon^2 \right).
\end{equation}
where $\mathbf{w}_1^{(i)}(n)= \GG(n) \PP_i \uu$, and where $\mathbf{G}(n)$ is defined in Eq.~\eqref{eq:def:G}. These transient states are plotted for $\mathbf{P}_1$ and $\mathbf{P}_2$, respectively, in Fig.~\ref{fig:lanford}(B) and (C). As $n$ tends to infinity the asymptotic linear response converges to $\mathbf{v}_{1}^{(i)}=\GG\PP_i \uu$, for $i=1,2,3$. In this particular case study, after three iterations of the transient linear response we observed convergence and the comparison to the invariant vector computed from Eq.~\eqref{eq:eig_prob} is shown as well in blue and red, in Fig.~\ref{fig:lanford}(B) and (C), respectively, were we observe that the linear response correction provides a faithful representation of the perturbed state.

\subsection{A double well potential with rotation}

\begin{figure}[H]
\centering
\includegraphics[scale=0.29]{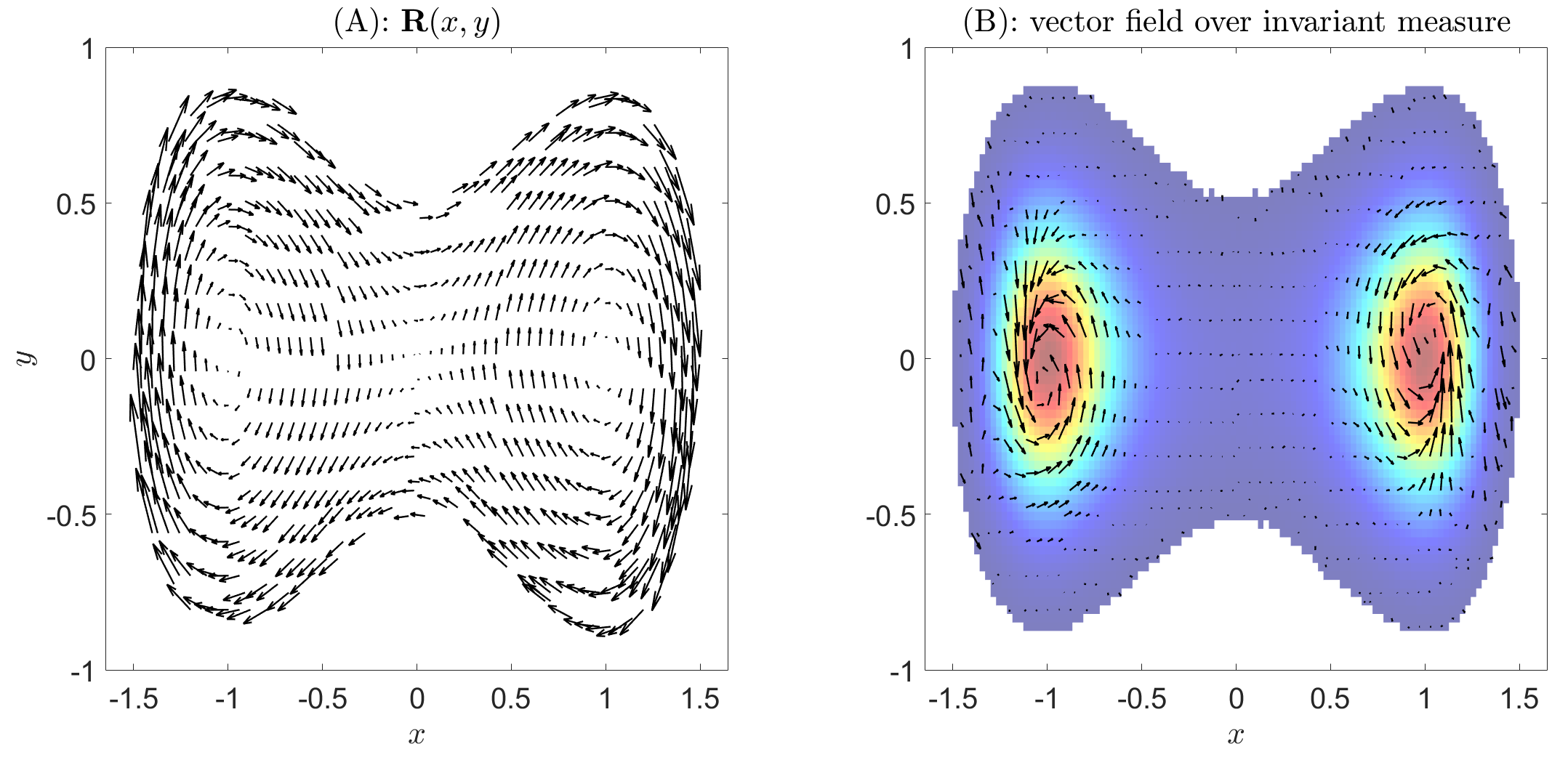}
\caption{\label{fig:dwp_rot1}\textbf{Double well potential with rotation: vector fields.} Panel (A) shows the applied rotational for $\mathbf{R}$ shown in Eq.~\eqref{eq:rotation}. Panel (B) shows the reconstructed perturbation vector field--- associated with $\mathbf{P}_s$--- superimposed to the unperturbed invariant measure.}
\end{figure}

The double well system considered above is a classical example of reversible dynamics satisfying detailed balance conditions, and thus with a vanishing entropy production. In this section, we add a rotational non-conservative force that drives the system out of equilibrium and yields a non vanishing entropy production. In particular, we here consider a rotational force 
$\mathbf{R}$ that does not alter the invariant density of the system but changes its dynamical properties by making it irreversible. The Fokker-Planck equation for the probability density $\rho$ associated to \eqref{eq:dwp_symm} can be written as $\partial_t \rho = -\nabla \cdot \mathbf{J}(\rho)$, where $\mathbf{J}(\rho) = \nabla \cdot( \mathbf{F\rho})-\frac{\sigma^2}{2}\nabla^2\rho$ represents the probability flux. In the double well system considered above, $\mathbf{J}(\rho_0)=0$ as the system is in detailed balance condition \cite{pavliotisbook2014}. Here, we introduce a non-conservative component $\mathbf{R}$ such that $\mathbf{J}(\rho_0) \neq 0$ but does not alter its invariant density $\rho_0$. Such a non-conservative component must satisfy $\nabla \cdot \mathbf{J}(\rho_0) =0$ \cite{pavliotisbook2014}:
\begin{equation}
    \nabla \cdot \left( \mathbf{R}(x,y)\rho_0(x,y) \right) = 0.
\end{equation}
There are many ways to choose the rotation $\mathbf{R}$, we take it as:
\begin{equation}\label{eq:rotation}
    \mathbf{R}(x,y) =\frac{1}{2} \binom{y}{x-x^3},
\end{equation}
for every $(x,y)$ in $\mathbb{R}^2$ and we represent it in Fig.~\ref{fig:dwp_rot1}(A) using a quiver-plot. 
The resulting double well system with rotation is \begin{equation}
    \dd \mathbf{x} = \left(\mathbf{F}\left( \mathbf{x}\right) + \mathbf{R}(\mathbf{x})\right)\, \dd t + \sigma\dd\mathbf{W}_t,
\end{equation}
where the noise strength is $\sigma=0.4$. We integrate the above equations for $10^7$ time units with a timestep of $10^{-2}$ . This provides a time series $\{(x_k,y_k) \}_{k=1}^{K}$ with $K = 10^8$ data-points starting at the origin, and recording snapshots every $0.1$ time units. As spatial domain we consider $[-2,2]\times [-1.5,1.5]$ and it is subdivided in the same way as in the previous example. A transition Markov matrix is obtained for a transition time of $\tau  = 0.1$ time units. 

Since the clockwise rotation $\mathbf{R}$ is the only process driving the system out of equilibrium to a state with a non vanishing entropy production, we expect that the entropy-production-minizing optimal perturbation will try to go against the rotational effect. We test this using the algorithm described in \cref{subsec:Entropy_prod} and then reconstructing the drift via the argument in \cref{subsec:drift_reconstr}. The optimal perturbation $\mathbf{P}_s$ is constructed with the explicit formula \eqref{eq:lagrange_matrix_formula} derived. The vector field associated with $\mathbf{P}_s$ is shown in Fig.~\ref{fig:dwp_rot1}(B) using a quiver-plot superimposed on top of the unperturbed invariant measure. A value of $\varepsilon=10^{-3}$ was employed. We observe that there are two main rotational features around the modes of the distribution, just as in Panel (A). Furthermore, we note that these rotational features gyrate counter-clockwise, opposite, as expected, to the applied rotation $\mathbf{R}$. The \textit{globally optimal} perturbation that minimizes the entropy production would be the one associated to the vector field $-\mathbf{R}$. This (large) perturbation would eliminate the rotational effect everywhere in phase space, restoring the detailed balance conditions and the equilibrium properties of the system and thus making the entropy production vanish. However, in a linear response approach, only perturbations of limited amplitude are allowed. The optimal perturbation found by the algorithm instead concentrates and optimizes the counter-rotational effect near the minima of the two wells, thus trying to reduce the entropy production where the system spends most time.

\begin{figure}[H]
\centering
\includegraphics[scale=0.29]{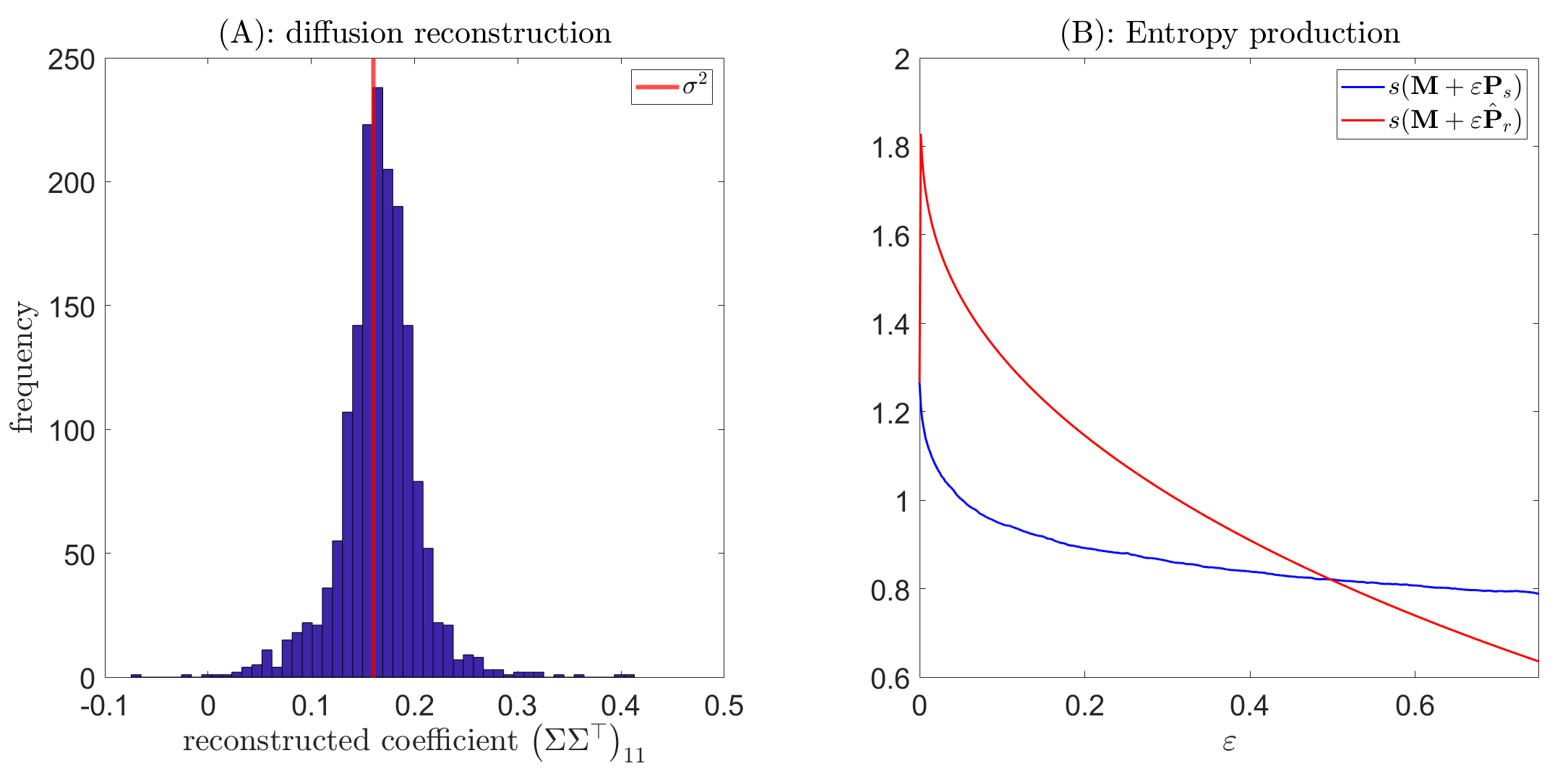}
\caption{\label{fig:dwp_rot2}\textbf{Double well potential with rotation: diffusion and entropy production.} Panel (A) shows a histogram of the values of the diagonal diffusion coefficient $\left(\Sigma\Sigma^{\top} \right)_{11}$ at each point in phase space, computed from the perturbed matrix $\mathbf{M}+\varepsilon\mathbf{P}$. The vertical red line indicates the value of the original unperturbed system. Panel (B) shows the entropy production associated with the matrices in the legend, for a range of values of $\varepsilon>0$.}
\end{figure}

To assess the diffusion introduced by the perturbation, we reconstructed the diffusion tensor $ \Sigma\Sigma^{\top}(x,y)$ for $\mathbf{M}+\varepsilon \mathbf{P}$. In Fig.~\ref{fig:dwp_rot2}(A), we are showing a histogram of the values of $\left( \Sigma\Sigma^{\top}\right)_{11}$ over all the domain. The modal value is located exactly at $\sigma^2$--- which is the analytical value of the unperturbed model---, suggesting that diffusion present in $\mathbf{M}$ is not affected by the introduction of $\mathbf{P}_s$.

It remains to check how the (linear) optimal perturbation  $\mathbf{P}_s$ differs from the perturbation matrix $\mathbf{P}_r$ defined in Eq.~\eqref{eq:reversibilized_perturbation} using the additive reversibilization of Markov matrices. First of all, we show in Appendix~\ref{app:additive_reversibilization_check} that the finite-difference drift reconstruction algorithm of Eq.~\eqref{eq:drift_recons} applied to $\mathbf{P}_r$ provides the \textit{globally optimal} rotational vector field $-\mathbf{R}$. We have verified that this perturbation fully restores equilibrium condition by numerically checking that $s(\mathbf{M}+\mathbf{P}_{s})=0$. However, this perturbation is large,  $\|\mathbf{P}_r\|_F\approx 4$, and thus far from a linear response regime. To compare the effects of the two perturbations for different various of $\varepsilon$, we introduced the normalized perturbation matrix $\hat{\mathbf{P}}_r = \mathbf{P}_r /\|\mathbf{P}_r\|_F $. For small values of $\varepsilon$, $\hat{\mathbf{P}}_r$ is not as good as $\mathbf{P}_s$ in reducing the entropy production. Indeed, this is shown in Fig.~\ref{fig:dwp_rot2}(B), where we plot the entropy productions $s(\mathbf{M}+\varepsilon \hat{\mathbf{P}}_{r})$ and $s(\mathbf{M}+\varepsilon\mathbf{P}_{s})$.

% \begin{figure}[H]
% \centering
% \includegraphics[scale=0.29]{images/revised_plot_JSP.png}
% \caption{\label{fig:dwp_rot}\textbf{A double well potential with rotation.} Panel (A): quiver plot of the vector field rotation $\mathbf{R}$ defined in Eq.~\eqref{eq:rotation}. Panel (B): vector field reconstruction associated with $\mathbf{P}_s$ superimposed to the invariant measure of the system. Panel (C): entropy production function of $\mathbf{M}+\varepsilon\mathbf{P}_s$ (blue) and $\mathbf{M}+\varepsilon \hat{\mathbf{P}}_r$ (red) vs. $\varepsilon$. }
% \end{figure}

\subsection{The stochastic Lorenz system}

The Lorenz 63 (L63) system is a low-dimensional truncation of the equations of atmospheric convection \cite{lorenzdeterministic1963}, whose dynamical properties have been extensively reported in the literature; we refer to the survey \cite{sparrow}. In its chaotic regime, while not displaying uniformly hyperbolic structure, it has been demonstrated on numerical grounds the existence of linear response, in the sense that small parameter forcings equate to proportional changes in the statistics \cite{reicklinear2002}, also under the perspective of Markov matrix \cite{Lucarini2016, SantosJSP}. The addition of stochastic noise into the L63 system models fast unresolved processes and has been analyzed in previous works, e.g. \cite{LUCARINI2026117540}. Here, we consider additive and isotropic noise yielding the following set of equations:
\begin{equation}\label{lorenz63}
\dot{\mathbf{x}}(t)=\mathbf{F}(\mathbf{x})=\begin{cases} s (y - x) + \sigma \dot{W}_{1,t} \\ x(r - z) - y + \sigma \dot{W}_{2,t} \\ xy -b z + \sigma \dot{W}_{3,t} \end{cases},
\end{equation}
where $\xx = (x,y,z)^\top$, $W_{i,t}$ denotes a one-dimensional Brownian motion. The model parameters are set to be $s=10$, $r=28$ and $b=8/3$, for which the system displays chaotic behavior \cite{Tucker2002}, and the noise intensity is set to $\sigma = 2$, following \cite{LUCARINI2026117540}.

To estimate the semigroup $\{ e^{t\LLL} \}_{t\geq 0}$, the system Eq.~\eqref{lorenz63} is integrated with a time step of $\mathrm{d}t = 5\cdot 10^{-3}$ time units  using an Euler-Maruyama method and initial condition $x=y=z=1$. After removing a transient of $100$ time units, snapshots where recorded every $0.05$ time units until a total of $K = 10^8$ snapshots were obtained. The trajectory $\{(x_k,y_k,z_k)\}_{k=0}^{K}$ is contained in the domain $ \mathcal{D} = [- 22,22]\times [-30,30]\times [-1,52]$. Such a domain is then divided into $N = 2^{18}$ equally sized, non-intersecting boxes $\{B_i\}_{i=1}^N$ of dimensions $0.69 \times 0.94 \times 0.88$. Non-populated boxes are discarded. Then, the Markov matrix $\mathbf{M}$ associated to \eqref{lorenz63} is computed using Eq.~\eqref{eq:transition_matrix} as before. Transition times values of $\tau=0.05,0.1$ and $0.2$ are investigated to examine the dependence of the results on such a parameter. We note that for all of the three values correlation functions can be well estimated up to 10 model units, as done previously in \cite{LUCARINI2026117540}. Since correlation functions are computed from the eigenpairs of the transition matrix, this suggests that the obtained matrices contain meaningful spectral information about the system \cite{chekroun2019c}.

\begin{figure}[H]
\centering
\includegraphics[scale=0.30]{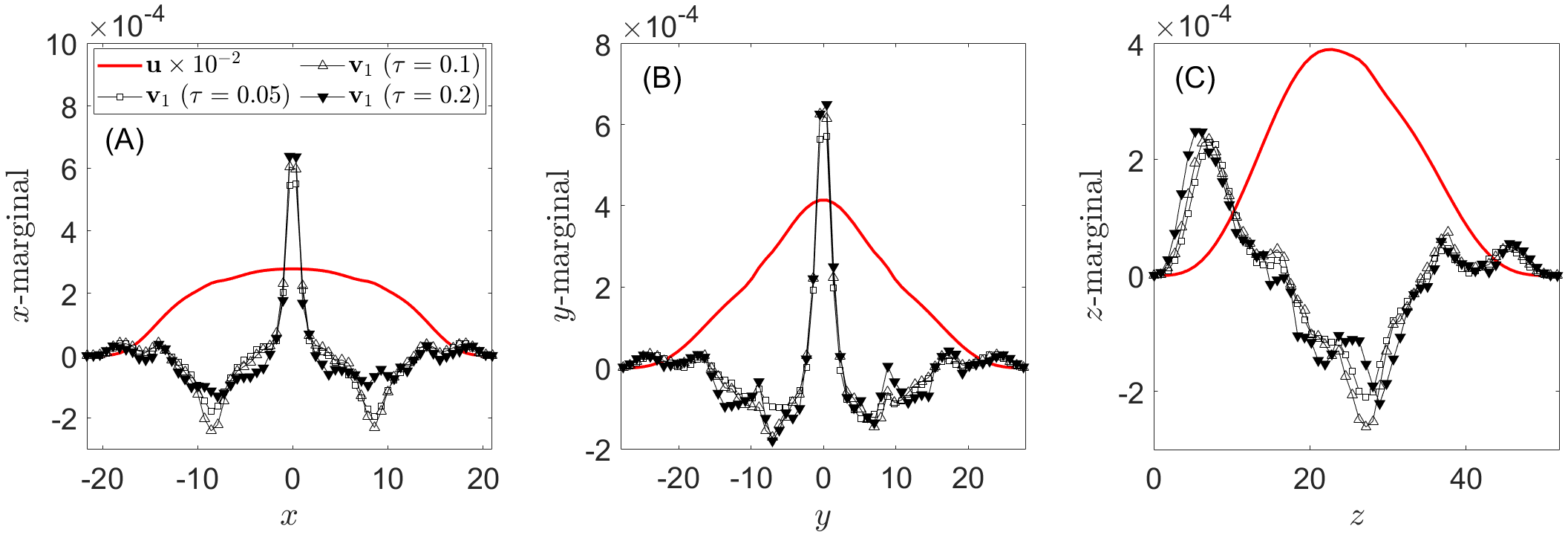}
\caption{\label{fig:l6311}\textbf{Linear response marginals for the stochastic L63 system.} Panels (A), (B) and (C) show the marginals in the $x$-, $y$- and $z$-directions, respectively, linear response of $\mathbf{M}$ as perturbed by $\mathbf{P}$. The different line formats correspond to difference values of transition time $\tau$, as indicated in the legend of panel (A). The red lines indicate the marginals of the unperturbed invariant measure and is scaled by $10^{-2}$ for visualization.}
\end{figure}

We consider the problem of finding the perturbation to the system that maximizes the entropy of the invariant distribution, following Sec.~\ref{subsec:entropy}. The formula~\eqref{eq:antown_expect} was applied to the transition matrices obtained for $\tau=0.05,0.1$ and $0.2$ to obtain their respective optimizing matrices together with their linear responses. In Fig.~\ref{fig:l6311} we show the marginal of the linear response for each of the system's variables: $x$-, $y$- and $z$-component in panels (A), (B) and (C), respectively. The red lines in each panel correspond to the marginals of the invariant probability vector $\mathbf{u}$ of $\mathbf{M}$ for $\tau = 0.1$, which does not change if $\tau$ changes. The vector $\mathbf{u}$ was scaled by $10^{-2}$ to make the visualization qualitatively comparable to the linear response plots. We observe two things. First of all, that the linear response is independent of $\tau$. Secondly, to maximize the entropy of the system, probability needs to be transferred to the center of the domain, concentrating along the $x=0$ and $y=0$ axis. This is counterintuitive, since increasing entropy should equate to flattening the distribution. On the $z$-coordinate, however, probability is transferred into lower values.

\begin{figure}[H]
\centering
\includegraphics[scale=0.255]{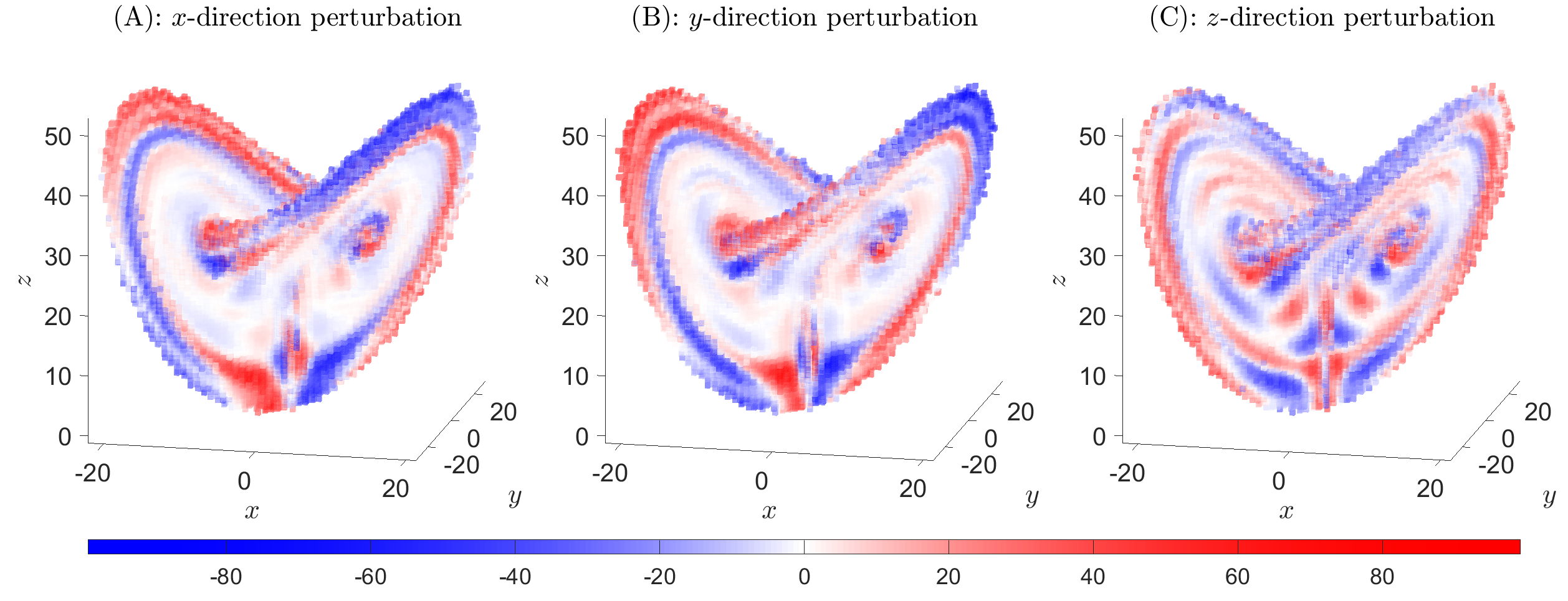}
\caption{\label{fig:l6312}\textbf{Perturbation maximizing entropy in the stochastic L63 system.} We show in panels (A), (B) and (C), the $x$-, $y$- and $z$-components of the perturbation vector field reconstructed from $\mathbf{P}$, where $\mathbf{P}$ maximizes the entropy of the system.}
\end{figure}

To examine how such a perturbation looks at the phase space level, we employ the scheme in Sect.~\ref{subsec:drift_reconstr} to reconstruct the vector field, for a perturbative value of $\varepsilon=10^{-3}$, which was seen to fall into the linear response regime. In Fig.~\ref{fig:l6312}, we show the $x$-, $y$- and $z$-coordinates of the perturbation vector field on the box covering of the domain (see panels (A), (B) and (C), respectively). For the $x$- and $y$-directions, the perturbation reconstruction is virtually symmetric. One can highlight the two red and blue features at the bottom of the attractor (low values of $z$), which could explain how the perturbation would transfer probability mass toward the lower part $x=y=0$ axis, as seen in the linear response plots of Fig.~\ref{fig:l6311} (panels (A) and (B)). In the upper part of the attractor (higher values of $z$) we also observe that the perturbation has intense red and blue colors. However, in comparison with the areas with low values of $z$, these are sparsely populated regions where the invariant vector $\mathbf{u}$ has low values.

\subsection{The projected stochastic L63 system}\label{sec:projected_l63}

In this section we consider the same stochastic L63 system as before, but we apply the methods to transition matrices computed on the $x$-$z$ projection of the timeseries data obtained in the previous section. The $x$-$z$ domain, $[-22,22]\times [-1,52]$, is then divided into $N = 2^{12}$ equally sized, non-intersecting boxes $\{B_i\}_{i=1}^N$ of dimensions $0.69 \times 0.88$. Then, the transition times of  $\tau=0.05,0.1$ and $0.2$ are selected to the Markov matrix associated to \eqref{lorenz63} is computed using Eq.~\eqref{eq:transition_matrix} as before. Note that because of the projection, we expect to lose the Markov property and, consequently, the spectral properties of the Markov matrices will depend on $\tau$. In this section we aim at illustrating this fact.

We obtained the perturbation matrix $\mathbf{P}$ such that the resulting invariant vector maximizes the KL-divergence with respect to the unperturbed vector $\mathbf{u}$; see Eq.~\eqref{eq:rel_entr}. Then, the perturbation size of $\varepsilon = 0.01$ is chosen so that $\MM + \varepsilon \PP$ is a (perturbed) Markov matrix with perturbed invariant measure $\mathbf{v}$ satisfying:
\begin{equation}\label{eq:eig_prob11}
    \left( \MM + \varepsilon \PP\right)\mathbf{v}=\mathbf{v}=\uu + \varepsilon \vv_1 + \mathcal{O}\left(\varepsilon^2 \right),
\end{equation}
where we recall that the (asymptotic) linear response is given by $\vv_1=\mathbf{G}\mathbf{P}\mathbf{u}$. We recall that the matrix $-\mathbf{P}$ would also maximize KL-divergence from $\mathbf{u}$, but we shall omit the analyses for the opposite-sign matrix.

\begin{figure}[H]
\centering
\includegraphics[scale=0.27]{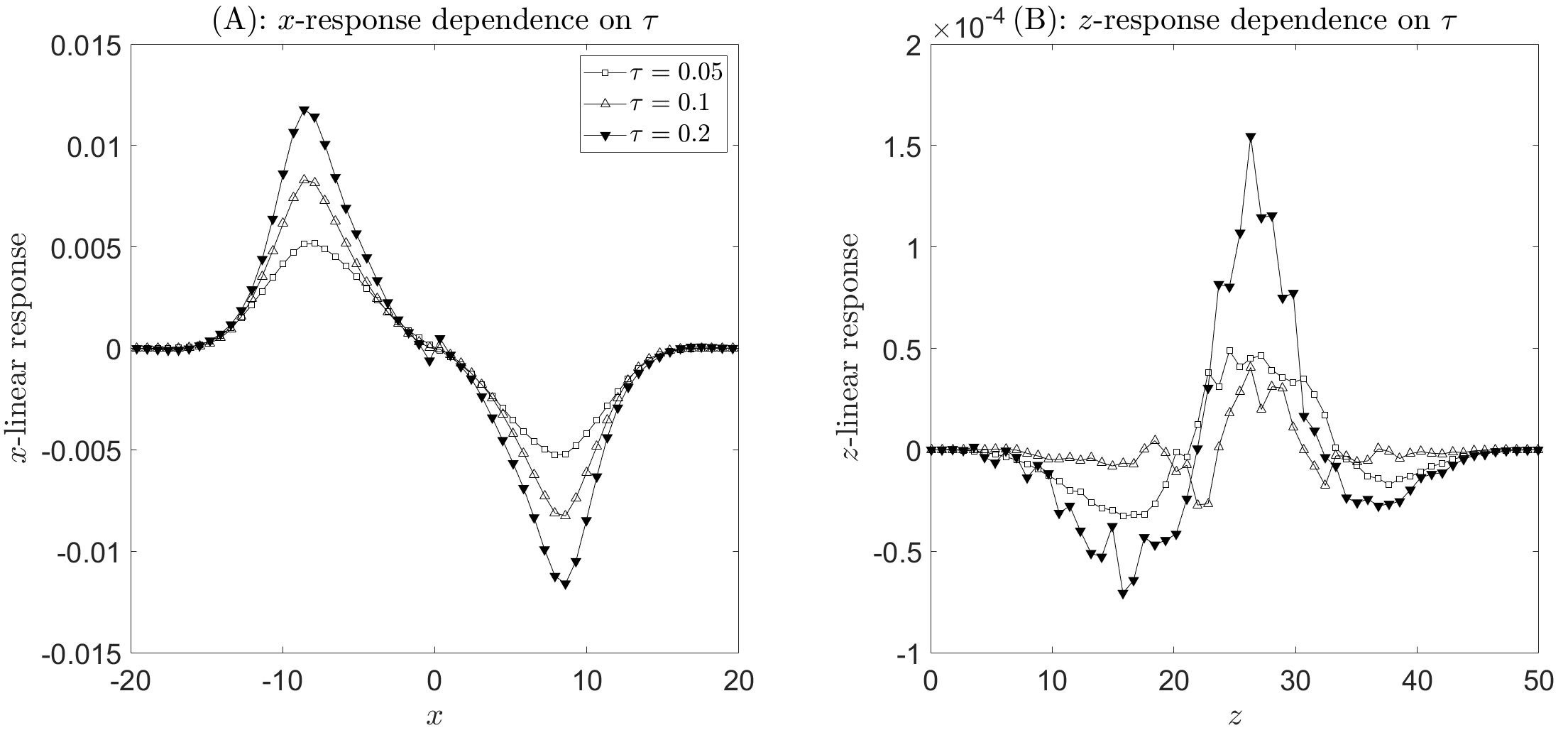}
\caption{\label{fig:l63_2D_11}\textbf{Dependence on $\tau$ for the projected L63 system.} Panels (A) and (B) show the marginals in the $x$- and $z$-directions, respectively, of the linear response of $\mathbf{M}$ as perturbed by $\mathbf{P}$. The different line formats correspond to different values of transition time $\tau$, as indicated in the legend of panel (A).}
\end{figure}

As mentioned earlier, because of the projected dynamics we expect that the Markov matrices constructed from a two-dimensional dataset will experience a non-trivial dependence on $\tau$ and, furthermore, the linear response to the perturbation matrices cannot be expected to remain invariant under $\tau$. Indeed, in Fig.~\ref{fig:l63_2D_11} we show in panels (A) and (B), the marginals of the scaled linear response along the $x$- and $z$-directions, respectively. We clearly observe--- as opposed to the full model in Fig.~\ref{fig:l6311}--- that the linear response severely depends on $\tau$.

\begin{figure}[H]
\centering
\includegraphics[scale=0.29]{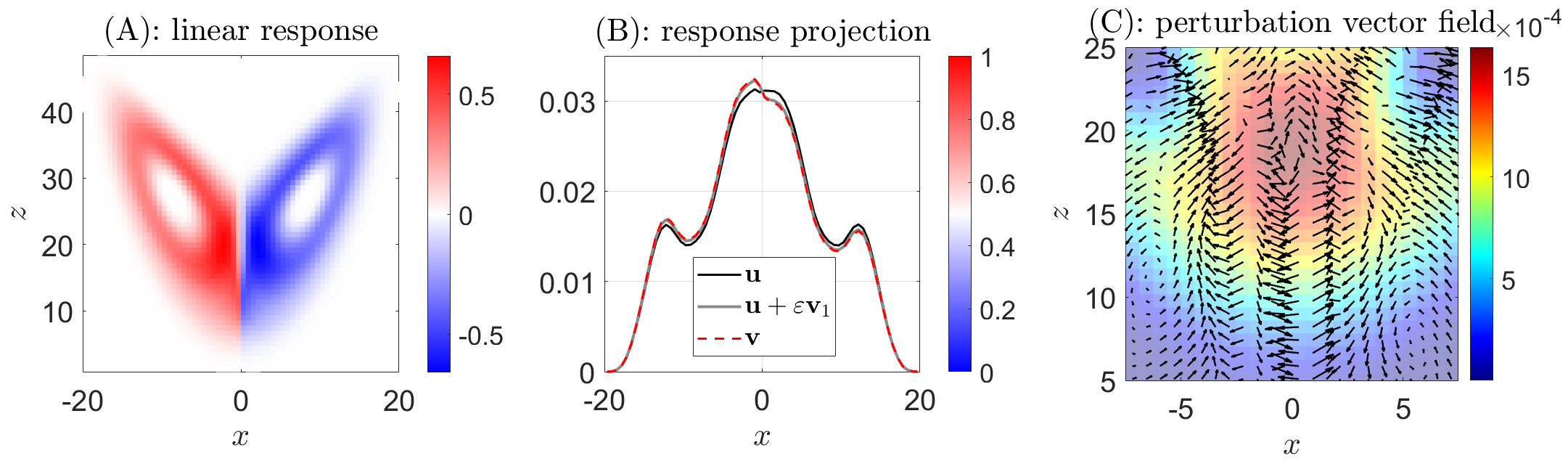}
\caption{\label{fig:l63_2D_1}\textbf{KL divergence in the projected L63 system.} Panel (A) shows the linear response $\mathbf{v}_1$ of $\mathbf{M}$ with respect to the perturbation $\mathbf{P}$ that maximizes KL divergence with respect to $\mathbf{u}$. Panel (B) shows the marginal response along the $x$-axis. The black line shows the unperturbed invariant vector, the red line shows the invariant vector of the matrix $\mathbf{M}+\varepsilon \mathbf{P}$ and the gray line shows the linear order correction to $\mathbf{u}$ due to the perturbation $\mathbf{P}$.}
\end{figure}

For illustration, we show in Fig.~\ref{fig:l63_2D_1}(A) the linear response vector $\mathbf{v}_1$ of the matrix $\mathbf{M}$ to the perturbation $\mathbf{P}$, for a transition time of $\tau=0.05$. Figure~\ref{fig:l63_2D_1}(B) contains the linear response marginal along the $x$-direction, revealing this time that the leading order correction of the invariant vector due to the perturbations is antisymmetric about the $x=0$ axis, in comparison to what was obtained in Fig.~\ref{fig:l6311}. The black and red lines show the unperturbed and perturbed invariant vectors, respectively. The gray line is the linear response correction, which coincides perfectly with the perturbed invariant vector, implying that we are within the linear response regime. This time, the perturbation tilts the weight of the distribution towards one of the lobes of the attractor.

% \begin{figure}[H]
% \centering
% \includegraphics[scale=0.29]{images/l63_ex2_t1.png}
% \caption{\label{fig:l631}\textbf{Entropy maximization in the L63 system.} Panel (A): perturbed invariant vector of the matrix $\mathbf{M}+\varepsilon\mathbf{P}_1$. Panel (B): linear response for the perturbation matrix $\mathbf{P}_1$. Panel (C): marginal distributions for the unperturbed and perturbed invariant vectors. The black line shows the invariant vector of $\mathbf{M}$ and the red shows the invariant vector of $\mathbf{M}+\varepsilon\mathbf{P}_1$. The grey line shows the linear order correction from Eq.~\eqref{eq:eig_prob11}.}
% \end{figure}

% \begin{figure}[H]
% \centering
% \includegraphics[scale=0.29]{images/l63_ex1_t1.png}
% \caption{\label{fig:l632}\textbf{KL-divergence maximization in the L63 system.} Panel (A): unperturbed invariant vector of the matrix $\mathbf{M}$. Panel (B): linear response for the perturbation matrix $\mathbf{P}_2$. Panel (C): marginal distributions for the unperturbed and perturbed invariant vectors. The black line shows the invariant vector of $\mathbf{M}$ and the red shows the invariant vector of $\mathbf{M}+\varepsilon\mathbf{P}_2$. The grey line shows the linear order correction from Eq.~\eqref{eq:eig_prob11}.}
% \end{figure}

\subsection{Unstable periodic orbit description of turbulent flows}\label{subsec:upo}
The role of invariant periodic solutions in shaping dynamical properties of chaotic and turbulent systems has been the subject of intense research \cite{Cvitanovic1988,ChaosBook}. On the one hand, the subcritical transition to turbulence can be recast in terms of a sequence of bifurcations of a non trivial originating invariant solution of the system \cite{Paranjape2023,Avila2013}. On the other hand, the chaotic dynamical behavior of turbulent systems can be described in terms of the underlying unstable periodic orbits which coexist with the chaotic attractor \cite{Kazantsev1998,KAWAHARA_KIDA_2001,MAIOCCHI2024}. Of interest in these findings is the possibility of providing a low order, effective stochastic description of chaotic systems in terms of (a subset of) the underlying unstable periodic orbits (UPOs) \cite{Maiocchi2022,Yalniz2021}. Here, as a last numerical example, we investigate the response properties of the reduced Markov Chain description of a turbulent flow provided in \cite{Yalniz2021}. We consider a $3$D Kolmogorov flow described by the following forced incompressible Navier-Stokes equations 
\begin{equation}
\label{eq: Kolmogorov flow}
    \partial_t \boldsymbol{v} + \boldsymbol{v} \cdot \nabla \boldsymbol{v} = - \nabla p + \nu \nabla^2 \boldsymbol{v} + \mathbf{F},
\end{equation}
on a box domain $\mathcal{D} =[0,L_x] \times [0,L_y] \times [0,L_z]$ with periodic boundary conditions. Here $\mathbf{x}= (x,y,z)\in \mathcal{D}$, $\boldsymbol{v}= (v_x,v_y,v_z) =\boldsymbol{v}(\mathbf{x}) $ is the velocity field of the flow, $p=p(\mathbf{x})$ is the pressure field and $\nu$ is the kinematic viscosity. Furthermore, the flow is forced with $\mathbf{F}(\mathbf{x})=  \Omega \sin \left(\frac{2\pi y}{L_y} \right)\mathbf{e}_x$ where $\mathbf{e}_x$ is the canonical vector in the $x$ direction and $\Omega$ is the amplitude of the perturbation. We set the parameters for the flow as $\nu =0.05$, $\Omega=1$, $L_x=L_y=2\pi$ and $L_z = \pi$.
\begin{figure}
     \centering
     \begin{subfigure}[b]{0.49\textwidth}
    \includegraphics[width=\textwidth]{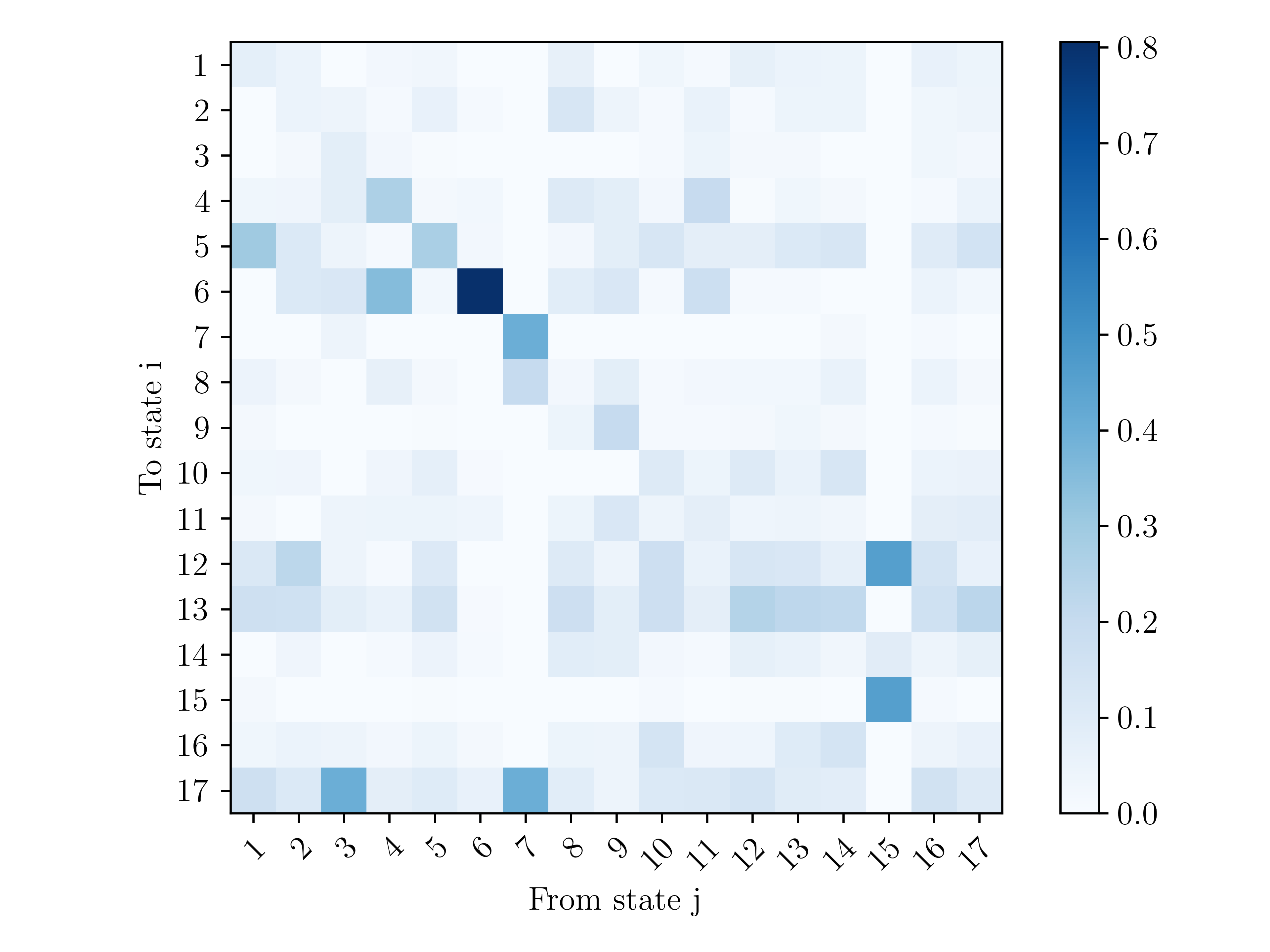}
    \caption{}
     \end{subfigure}
     \hfill
     \begin{subfigure}[b]{0.49\textwidth}
         \includegraphics[width=\textwidth]{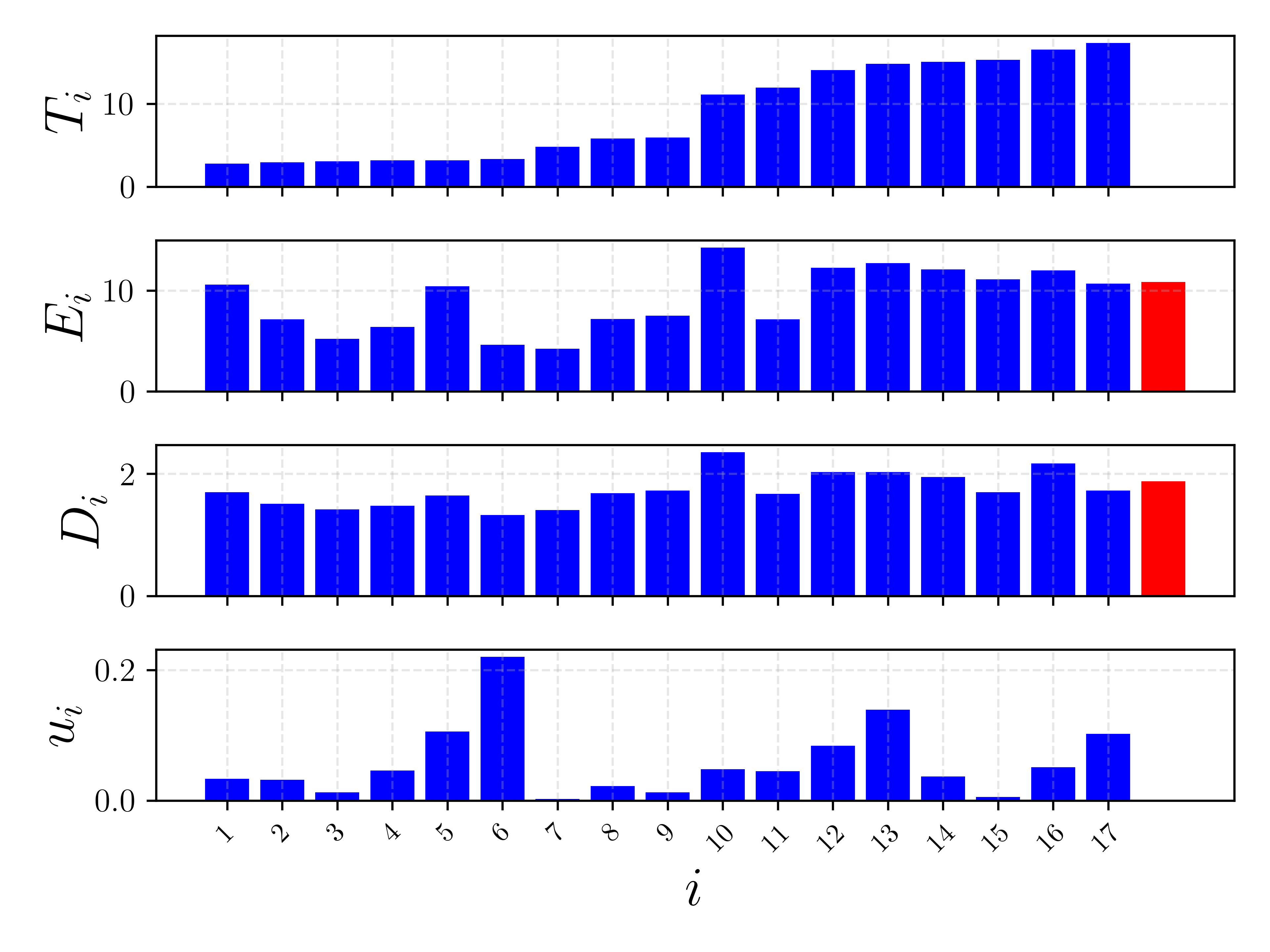}
         \caption{}
     \end{subfigure}
     \caption{Markov Chain description of the turbulent flow. Panel (a): Transition matrix $\mathbf{M}$ between periodic orbits. Panel (b): Properties of the unperturbed system. From top to bottom: Periods $T_i$, Energy $E_i$, Dissipation $D_i$ of the periodic orbits. Red bars correspond to the average \eqref{eq: weighted average UPO}. The bottom panel represents the unperturbed invariant density $\mathbf{u}$ of the Markov chain.}
     \label{fig: Unperturbed UPO system}
 \end{figure}
The above partial differential equation admits a laminar solution $\boldsymbol{v}_L = (\Omega \nu^{-1} (L_y/2\pi)^2 \sin(2\pi y/L_y),0,0)$ which is linearly stable for all values of $\nu$. However, one commonly observes transient (but long lasting) turbulent behavior if the system is triggered with finite amplitude perturbations. The authors in \cite{Yalniz2021} have shown that it is possible to obtain a reduced order model of the high dimensional turbulent dynamic in terms of the underlying unstable periodic orbits of the system. In particular, they managed to identify $17$ unstable periodic orbits $\gamma_i : \mathbb{R} \to \mathbb{R}^3$ such that $\gamma_i(t+T_i) = \gamma_i(t)$ $ \forall t$ where $T_i$ is the period of the $i-th$ unstable periodic orbit. The shadowing properties of the turbulent trajectory, i.e., the amount of time that the turbulent trajectory spends in the neighborhood of an unstable orbit, have been assessed with topological data analysis techniques by comparing the difference between the shapes of the periodic orbits and the turbulent trajectory.  
The resulting Markov Chain description corresponds with the identification of the periodic orbits as fundamental states of the system and with modeling the turbulent dynamic as a stochastic hopping between them. The statistics of such hoppings are given in terms of a transition matrix $\mathbf{M}=(M_{ij})$ describing the transition rates from the $j-th$ periodic orbit to the $i-th$ periodic orbit. The transition matrix is obtained by carefully studying the shadowing properties of the turbulent trajectory, see details in \cite{Yalniz2021}, and is reported here in Fig.~\ref{fig: Unperturbed UPO system}. 

The validity of this reduced order model is corroborated by the fact that statistics of the turbulent flow can be reconstructed from the Markov Chain description. In particular, given any observable of the system $h(\mathbf{x})$, such as the kinetic energy $E(\mathbf{x}) = |\boldsymbol{v}(\mathbf{x})|^2 / 2$ or the dissipation $D(\mathbf{x})= \nu |\nabla \times \boldsymbol{v}|^2$, it is possible to associate to each periodic orbit an average value
\begin{equation}
    \bar{h}_i = \frac{1}{T_i}\int_0^{T_i} h(\boldsymbol{\gamma}_i(s)) \mathrm{d}s.
\end{equation}
This way we construct the coarse-grained observable vector $\mathbf{h}=(\bar{h}_i)$ in $\mathbb{R}^{17}$. The invariant vector $\mathbf{u}=(u_i)$ of the coarse grained system satisfies $\mathbf{M}\mathbf{u}=\mathbf{u}$, and the averages of observables of the turbulent flow can be reconstructed from the Markov matrix picture by considering suitably weighted expectation values \cite{Yalniz2021}
\begin{equation}
\label{eq: weighted average UPO}
    \left\langle \mathbf{h} \right\rangle = \frac{\sum_{i=1}^{17}u_iT_i \bar{h}_i}{\sum_{j=1}^{17}u_jT_j} \coloneqq \sum_{i=1}^{17} w_i \bar{h}_i,
\end{equation}
where the weights $w_i = u_i T_i / (\sum_{j=1}^{17} u_j T_j)$ take into account the average residence time $\tau_i = u_i T_i$ of the turbulent trajectory on each periodic orbit. We highlight that Eq.~\eqref{eq: weighted average UPO} is a coarse grained approximation of the expectation value of the actual observable $h(\mathbf{x})$ evaluated on the turbulent trajectory.

We now consider a situation where a perturbation is applied to the system, such as the transitions between the unperturbed system's UPOs are described by a perturbed transition matrix $\mathbf{M}_\varepsilon = \mathbf{M} + \varepsilon \mathbf{P}$ where $\varepsilon$ is a small parameter. In the long time limit, the perturbed system will reach its statistical steady state characterized by the invariant distribution  $\mathbf{u}+\varepsilon \mathbf{v_1}$, with $\mathbf{v}_1$ being the asymptotic change to the invariant distributions given by \eqref{eq:assympt_linearresponse} as $\mathbf{v}_1 = \mathbf{GPu}$. Expectation values in the perturbed system can be written in a linear approximation as $\langle \mathbf{h} \rangle_\varepsilon = \langle \mathbf{h} \rangle + \varepsilon \sum_{i} \bar{h}_i \delta w_i$, where the perturbations to the weights is given by
\begin{equation}
\label{eq: perturbation weights}
    \delta w_i = \frac{T_i v_{1,i}}{\sum_{j=1}^{17} u_j T_j} - w_i \frac{\sum_{j=1}^{17} v_{1,j}T_j}{\sum_{j=1}^{17} u_jT_j} .
\end{equation}

We seek for the perturbation matrix $\mathbf{P}$ that maximizes the change in expectation value $\sum_i \bar{h}_i \delta w_i$. From the previous expression it is simple to see that the change in the weights can be written as $\delta w_i = ( \mathbf{A} \mathbf{v}_1)_i$, where
the matrix $\mathbf{A}$ is
\begin{equation}
    A_{ij} = \frac{T_j \delta_{ij} - w_iT_j}{\sum_{k=1}^{17} u_k T_k}
\end{equation}
and $\delta_{ij}$ is the Kronecker delta. We are then trying to find the optimal perturbation that solves the optimization problem
\begin{equation}
\operatorname{arg max}_{\mathbf{P}}\sum_i \bar{h}_i \delta w_i =  \operatorname{arg max}_{\mathbf{P}}\sum_{i,j} \bar{h}_i A_{ij} v_{1,j} = \operatorname{arg max}_{\mathbf{P}} \mathbf{f}^\top (\mathbf{GPu})
\end{equation}
where we have defined $\mathbf{f} = \mathbf{A}^\top \mathbf{h}$.
 This shows that this optimization problem falls into the same class given by ~\eqref{eq:max_entropy}, and we can then apply the same formula ~\eqref{eq:antown_expect} for the optimal perturbation $\mathbf{P}$.
\begin{figure}
     \centering
     \begin{subfigure}[b]{0.49\textwidth}
    \includegraphics[width=\textwidth]{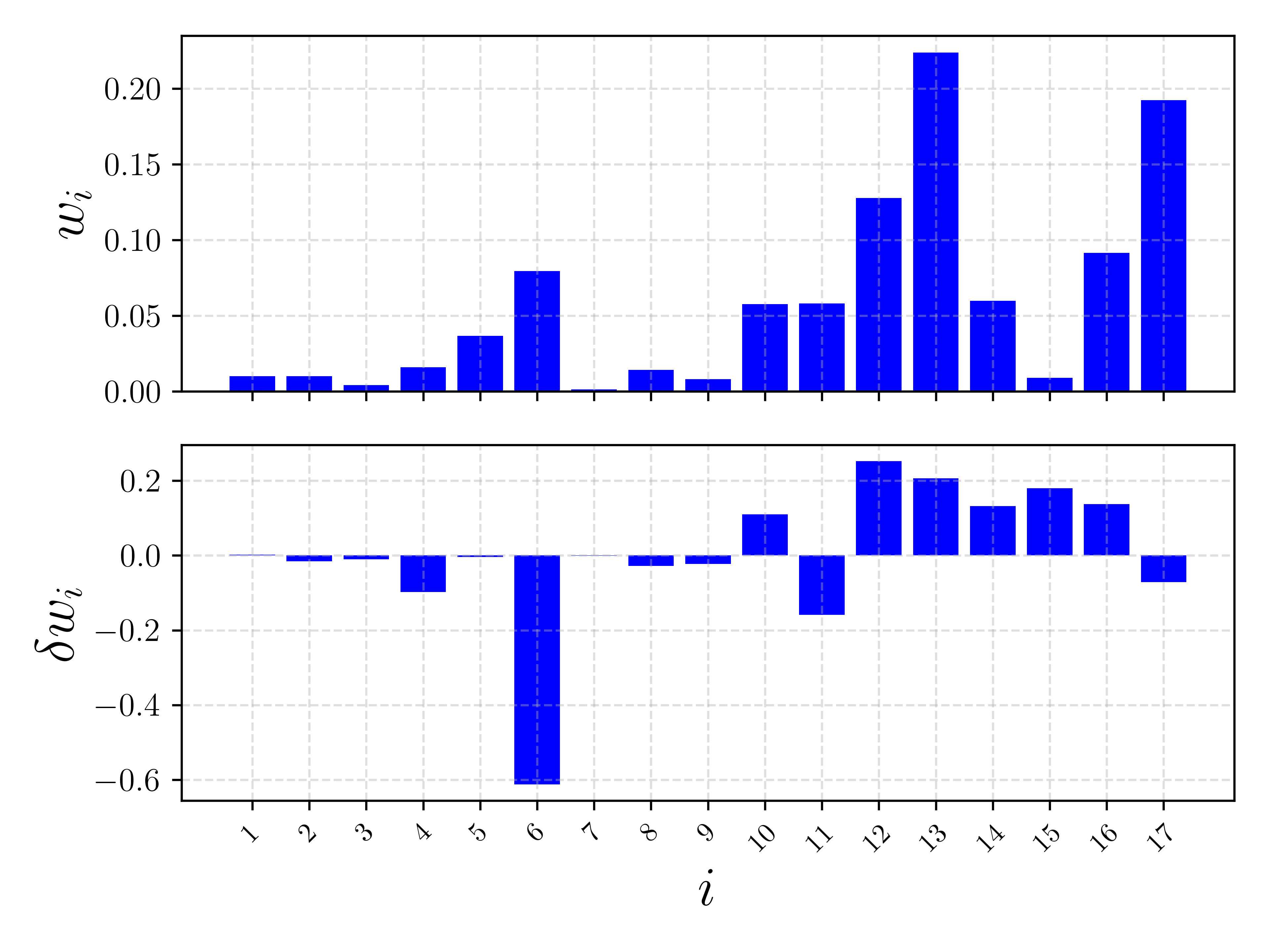}
    \caption{}
     \end{subfigure}
     \hfill
     \begin{subfigure}[b]{0.49\textwidth}
         \includegraphics[width=\textwidth]{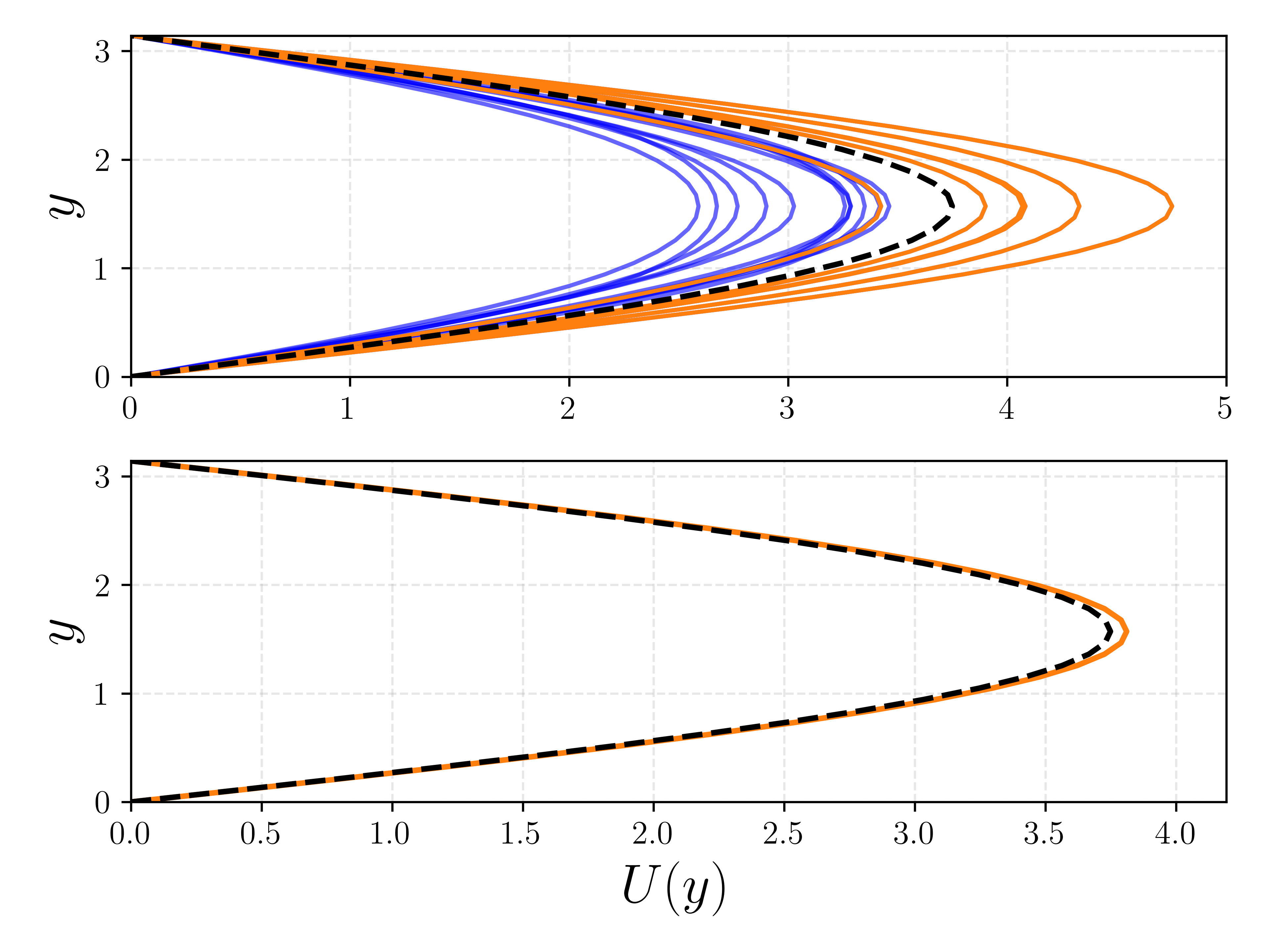}
         \caption{}
     \end{subfigure}
     \caption{Properties of the perturbed system. Panel (a): weights of the system. Top: weights of the unperturbed system. Bottom: change in the weights due to the perturbation. Panel (b): Velocity profiles $U(y)$ of the system. Top: continuous lines correspond to the UPOs profiles, the dashed line to the average profile of the turbulent trajectory obtained with \eqref{eq: weighted average UPO}. We have highlighted in orange the velocity profiles of the periodic orbits selected by the optimization algorithm as explained in the main text. Bottom: Effect of the perturbation on the average profile. The dashed (solid) line corresponds to the unperturbed (perturbed) average velocity profile. Here $\varepsilon = 0.05$ and profile are shown only for $y \in [0,\pi]$ because of symmetry properties of $U(y)$.}
     \label{fig: Perturbed UPO system}
 \end{figure}
 \\ \\
Here, we consider the optimal perturbation $\mathbf{P}$ that maximizes the average kinetic energy of the system $E(\mathbf{x})$ evaluated on the turbulent trajectory. We remark that, in a linear response approach, the perturbation that instead minimizes the energy is simply $- \mathbf{P}$. Panel (a) in Fig.~\ref{fig: Perturbed UPO system} shows the result of such perturbation. Here, we have checked that linearity, and thus the validity of \eqref{eq: perturbation weights}, holds for $\varepsilon \in (0,0.05]$. The perturbation strongly penalizes the weight $w_6$ which corresponds to a state with a high occupancy ($u_6$ corresponds to the maximum of the invariant density) but with quite low energy compared to the other states. For a similar reason, state $11$ gets less populated than $10$ because, despite having $w_{10} \approx w_{11}$, the energy of the two states differs substantially. In general, the perturbation favors states with a longer period, which correspond to by and large high energy states. The perturbation acts in a non trivial way, as, for example, mass leaks from state $17$, corresponding to the longest and quite energetic periodic orbit (in fact with the second largest weight), to favor transitions to other states in order to maximize the energy of the system. 

We remark that, as opposed to the previous section where an Ulam approach has been taken, this analysis does not allow for an immediate reconstruction of the optimal perturbation field to be applied to \eqref{eq: Kolmogorov flow}. This is coherent with the fact that, within an Ulam approach, we seek for an approximation of the semigroup of the dynamics using a finite basis of indicator functions $\chi_{B_i}$ where $B_i$ are the \textit{cells} that we have used for partitioning the phase space, see section \ref{sec:ulam}. Instead, in this section we consider a coarse-graining of the turbulent dynamics in terms of unstable periodic orbits, which are global quantities with a non trivial geometry in space\footnote{We refer the reader to the supplementary video in \cite{Yalniz2021} for a visual representation of the dynamics of the turbulent trajectory and the shadowing of the unstable periodic orbits.}. Below, we show that the optimal perturbation $\mathbf{P}$ identified above pushes the system towards a bundle of periodic orbits with similar spatial properties. To characterise the spatial properties of the UPOs, we follow \cite{Yalniz2021} and we concentrate on the spatially dependent velocity profile of solutions of \eqref{eq: Kolmogorov flow}
\begin{equation}
    U(y) = \frac{1}{L_x L_z}\int v_x(x,y,z) \mathrm{d}x\mathrm{d}z.
\end{equation} Panel (b, top) of \cref{fig: Perturbed UPO system} shows the average velocity profiles of the UPOs as continuous lines. As previously noted in \cite{Yalniz2021}, we observe that the velocity profiles of the UPOs can be quite different. In dashed line, instead, is the average velocity profile of the turbulent trajectory obtained using the Markov chain approximation \eqref{eq: weighted average UPO}. To investigate the effect of the perturbation $\mathbf{P}$ on the average velocity profile, we have highlighted in orange the UPOs that are favored by the perturbation, that is, the orbits corresponding to a positive weight change $\delta w_i >0$. We observe that the subset of UPOs selected by the perturbation correspond to velocity profiles that are higher than average. There is one exception to the above: the periodic orbit $i=15$ has a velocity profile below the average one . However, the weight of this state is very small due to a very low $u_{15} \approx 0.005$ invariant density. The overall effect of the perturbation on the average velocity profile --- see panel (b,bottom) of \cref{fig: Perturbed UPO system} --- correspond to stretching the peak of the profile around $y=\pi/2$.

\section{Conclusion}\label{sec:conclusion}

Discrete Markov chains are a natural tool for the statistical and ergodic analysis of dynamical systems and stochastic differential equations. Furthermore, the linear response theory of statistical physics is straightforwardly formulated in linear algebraic terms and serves to design approximation and optimization methods \cite{Lucarini2016,Antown2018}. In this work we provided linear optimization algorithms so that, given a Markov matrix, entropy, KL divergence and entropy production are minimized or maximized by adding a suitable perturbation. Moreover, a numerical link was presented between such perturbations at the matrix level and perturbations at the vector field level. We believe that this work, together with earlier ones, provide a basis for future research.

The topology of the network of states of a Markov chain and their connectivity heavily influence the location of its eigenvalues \cite{bremaud_markov_chains}. A more general set of algorithms should be constructed on top of the constraints (C1) to (C4)--- see Assumption~\ref{assum:perturb} and below Eq.~\eqref{eq:reversibilized_perturbation}---, so that a specific sparsity mask is imposed. At a first glance, linear functionals like that in Eq.~\eqref{eq:entropy_production} should be the easiest framework on which a specified transition network is imposed. As it is now, the present theory does not allow for the creation of new transitions. While creating new transitions might be unphysical in the context of Ulam discretizations, it might be relevant in other kinds of Markov chains or networks where there is no strict notion of continuity. Relaxing the sparsity constraint on top of considering nonlinear optimization functionals should bring new insights.

Ulam's method revisited here in Sect.~\ref{sec:ulam} suffers from the curse of dimensionality, due to the exponential growth in boxes needed to cover the relevant attractors or domains. This makes the method intractable for the numerical discretization of high-dimensional Fokker-Planck operators. Transition matrices can nevertheless be defined in high-dimensional settings by constructing kernel-smoothed transfer or Koopman operators \cite{Froyland2021}. Under this setting, it is possible to compute optimal perturbations and retrieve their associated perturbation vector fields \cite{froyland2026datadrivenmethodscomputationoptimal}. While Markov matrices in projected phase-spaces still provide useful information \cite{chekroun2019c,Tantet2015}, relating such matrices with an underlying flow is a hard task due to memory effects--- as was discussed in Sect.~\ref{sec:projected_l63} herein---. An effective drift reconstruction following Sect.~\ref{subsec:drift_reconstr} could be applied to design reduced order models from a Markov matrix perspective. Indeed, from partial observations one can estimate state transition rates following Ulam's algorithm--- shown here in Sect.~\ref{sec:ulam}--- and compute an effective vector field. Once a vector field is obtained, one could assess to which extent it can be used for prediction and, hence, quantify the role of the unresolved variables.

\paragraph{Acknowledgements}\quad The authors thank Gökhan Yalnız for sharing the data regarding the unstable periodic orbits for the Kolmogorov flow.  MSG acknowledges the partial support provided by the Horizon Europe project Past2Future (Grant No. 101184070). N.Z. has been supported by the Wallenberg Initiative on Networks and Quantum Information (WINQ) and by the European Space Agency Project PREDICT (Contract 4000146344/24/I-LR). GC is a member of GNAMPA of the Italian Istituto Nazionale di Alta Matematica (INdAM).

\pagebreak

\pagebreak
\appendix

\section{Finite-difference reconstruction}\label{app:finite_difference}

Here we shall justify Eq.~\eqref{eq:drift_recons} for the reconstruction of perturbation vector fields from perturbation matrices $\mathbf{P}$. The first assumption is that the matrix-logarithm of $\mathbf{M}+\varepsilon\mathbf{P}$ approximates the generator of the perturbed process. The difference with the matrix-logarithm of $\mathbf{M}$ should then give the approximation of the perturbation matrix $\mathbf{Q}$ such that:
\begin{equation}\label{eq:exp_derivation}
\tau \left[\left(\widehat{\mathbf{F}}_\varepsilon - \widehat{\mathbf{F}}\right)(\mathbf{c}_i) \right]_k \approx \sum_{j=1}^N \left[\log\left( \mathbf{M}+\varepsilon\mathbf{P} \right)- \log \left(\mathbf{M} \right)\right]_j \left(c^{(k)}_j-c^{(k)}_i\right)\approx \varepsilon \sum_{j=1}^N Q_{ji} \left(c^{(k)}_j-c^{(k)}_i\right) +\mathcal{O}(\varepsilon\tau),
\end{equation}
for each $k=1,\ldots,d$. We have to show that $\mathbf{Q} = \mathbf{P}$.

The second assumption is that the transition time $\tau$ is small enough so that the transition matrix $\mathbf{M}$ can be written as:
\begin{equation}
    \mathbf{M}  = \mathbf{I} + \tau \mathbf{H} + \mathcal{O}(\tau^2),
\end{equation}
for some matrix $\mathbf{H}$. This formula is to say that $\mathbf{M}$ is diagonally dominant and that, as $\tau$ decreases, it becomes closer to the identity matrix. The third assumption is that $\mathbf{M}$ is invertible. With all these three assumptions, the following perturbative relations hold:
\begin{subequations}
\begin{align}
    \mathbf{M}^{-1} &= \mathbf{I} - \tau \mathbf{H} + \mathcal{O}(\tau^2), \\
    \log \left( \mathbf{M}\right) &=\tau\mathbf{H}+\mathcal{O}(\tau^2).
\end{align}
\end{subequations}
Then we have:
\begin{subequations}
    \begin{align}
        \log \left( \mathbf{M} + \varepsilon\mathbf{P} \right) &= \log \left( \mathbf{M}\left(\mathbf{I} +\varepsilon\mathbf{M}^{-1}\mathbf{P}\right) \right) \\&= \log \left( \mathbf{M} \right) + \log \left(\mathbf{I} + \varepsilon \mathbf{M}^{-1}\mathbf{P} \right) \nonumber \\  &\quad+ \log \left( \mathbf{M}\right)\log \left( \mathbf{I+\varepsilon\mathbf{M}^{-1}\mathbf{P}} \right)-\log \left( \mathbf{I+\varepsilon\mathbf{M}^{-1}\mathbf{P}}\right)\log \left( \mathbf{M}\right) + \mathcal{O}(\varepsilon \tau) \\
        &=  \log\left(\mathbf{M} \right) + \varepsilon \mathbf{M}^{-1}\mathbf{P} + \varepsilon\tau \mathbf{H}\mathbf{M}^{-1}\mathbf{P} - \varepsilon\tau \mathbf{M}^{-1}\mathbf{P}\mathbf{H} + \mathcal{O}(\varepsilon \tau) \\ &= \log\left(\mathbf{M} \right)+ \varepsilon \left( \mathbf{I}-\tau \mathbf{H}\right)\mathbf{P} + \mathcal{O}(\varepsilon \tau)\\ &= \log\left(\mathbf{M} \right) + \varepsilon \mathbf{P} + \mathcal{O}(\varepsilon \tau).
    \end{align}
\end{subequations}
Finally, subtracting $\log\left(\mathbf{M} \right)$ and inserting into Eq.~\eqref{eq:exp_derivation} gives that $\mathbf{Q}=\mathbf{P}$.

\section{Robustness tests}
\subsection{Double-well potential reconstruction: resolution dependence}\label{app:resolution}
\begin{figure}[H]
\centering
\includegraphics[scale=0.29]{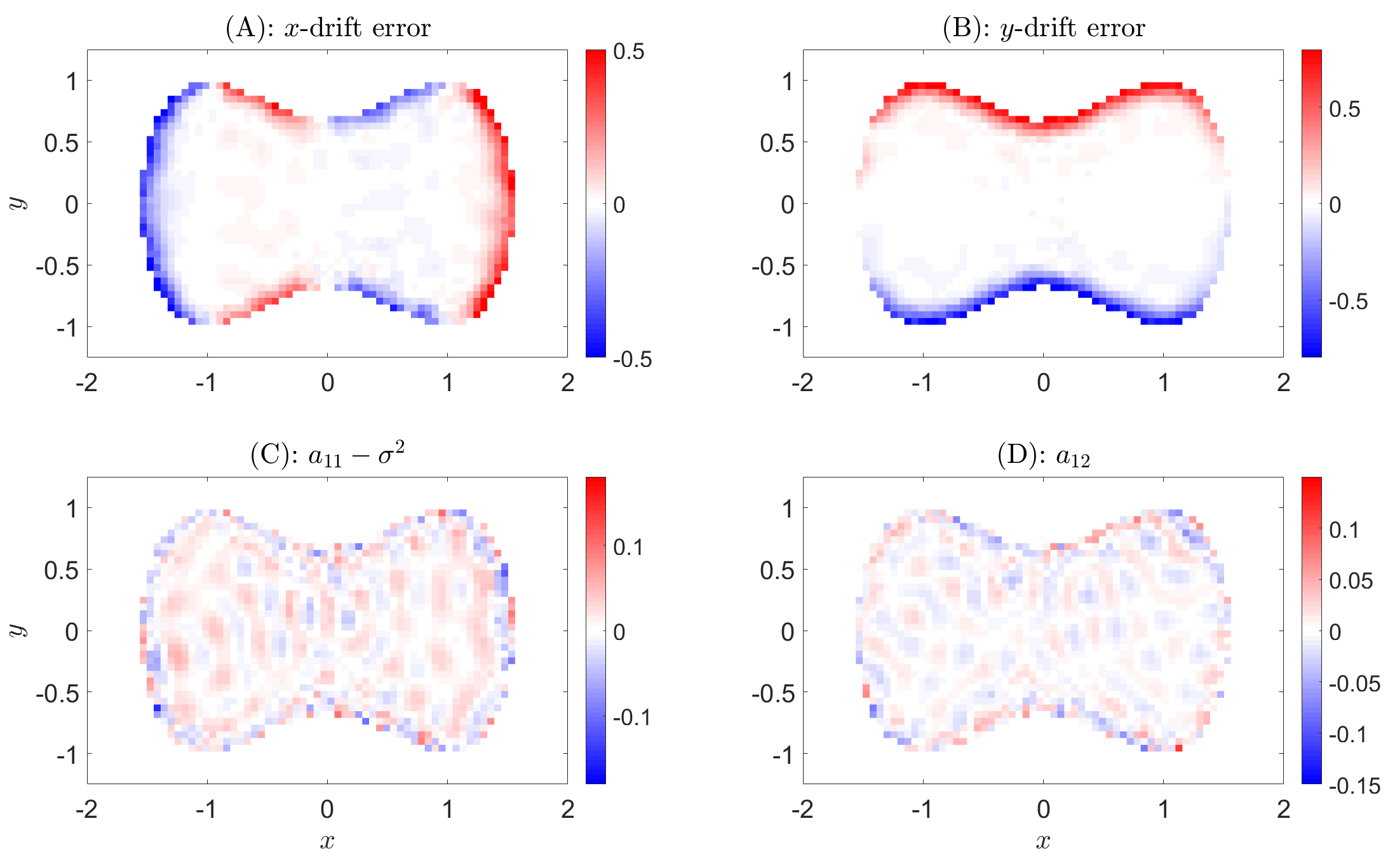}
\caption{\label{fig:dwp1_rec_depth12}\textbf{Double well dynamics reconstruction.} Panels (A) and (B) show the difference between the actual vector field and the reconstructed one for the $x$- and $y$-directions, respectively. On the bottom row, panel (C) shows the difference between the reconstructed and true first-diagonal component of the diffusion tensor, $\left(\Sigma\Sigma^\top \right)_{11}$. Panel (D) shows the value of the reconstructed off-diagonal entry of the diffusion matrix, $\left(\Sigma\Sigma^\top \right)_{12}$ (the true value vanishes everywhere). The resolution is of $2^{12}$ boxes covering the domain $[-2,2]\times [-1.5,1.5]$.}
\end{figure}

\begin{figure}[H]
\centering
\includegraphics[scale=0.29]{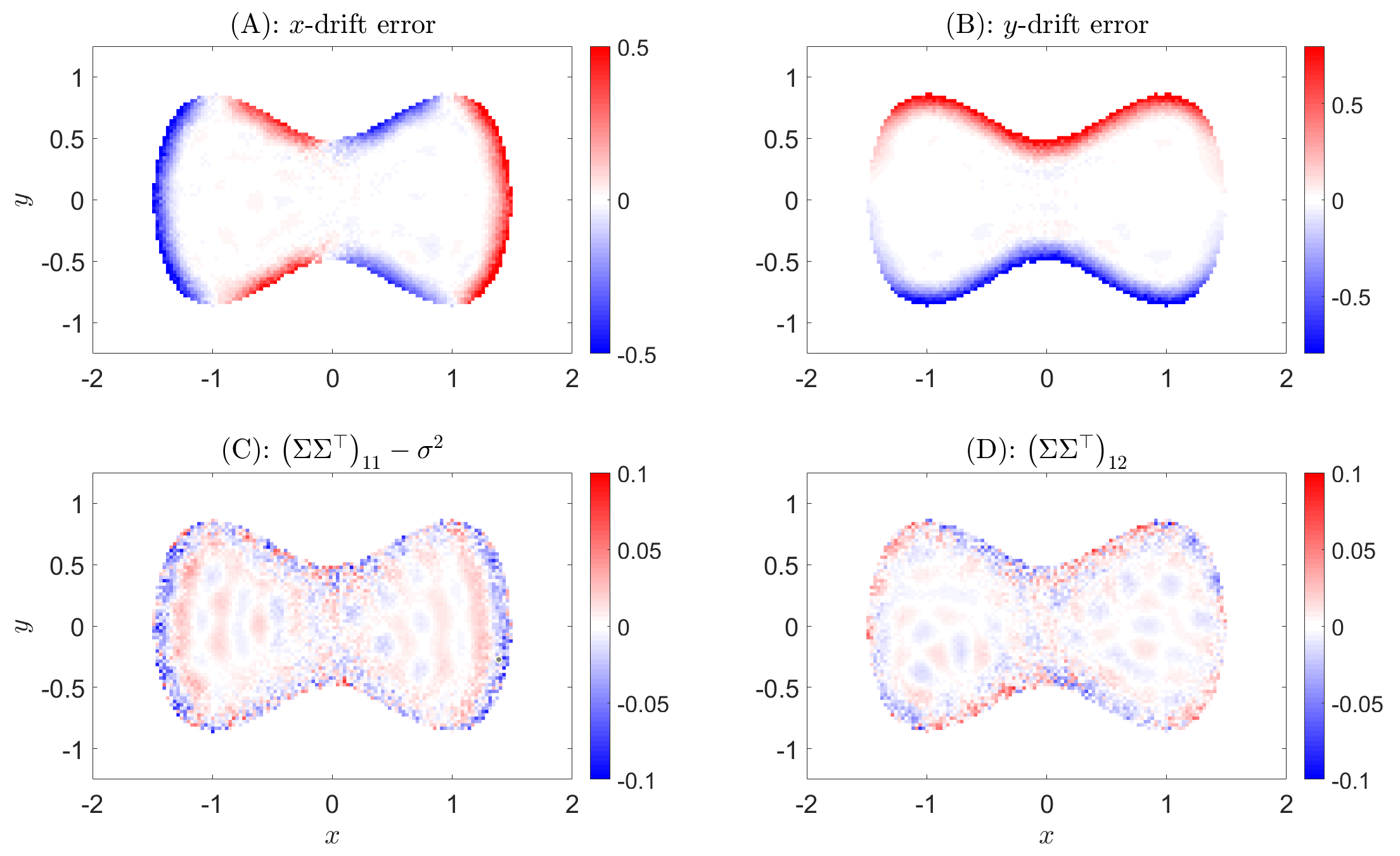}
\caption{\label{fig:dwp1_rec_depth14}\textbf{Double well dynamics reconstruction: low resolution.} Same as Fig.~\ref{fig:dwp1_rec_depth12}, but for a resolution is of $2^{14}$ boxes covering the domain $[-2,2]\times [-1.5,1.5]$.}
\end{figure}

In this appendix we examine the robustness of the reconstruction protocol to spatial resolution. In Fig.~\ref{fig:dwp1_rec_depth12} we show the error reconstruction of the perturbation vector field is shown in panels (A) and (B) and, in panels (C) and (D), we show the error in computing the diffusion tensor. Here we only plot the entries $\left( \Sigma \Sigma ^{\top}\right)_{11}$ and $\left( \Sigma \Sigma ^{\top}\right)_{12}$ of the diffusion tensor; panels (C) and (D), respectively. The resolution is that of $2^{12}$ boxes covering the domain $[-2,2]\times [-1.5,1.5]$. The same is shown in Fig.~\ref{fig:dwp1_rec_depth14}, but with a higher resolution of $2^{14}$ boxes.

\begin{figure}[H]
\centering
\includegraphics[scale=0.29]{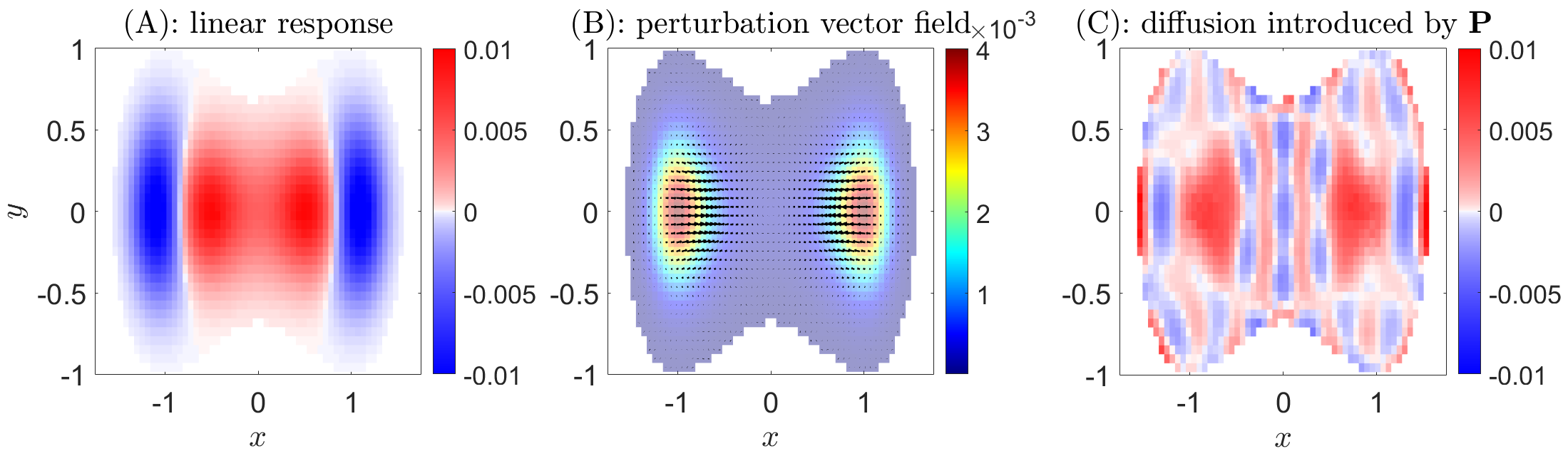}
\caption{\label{fig:dwp2_resp_lowres}\textbf{Expectation value maximization: perturbation and linear response.} Panel (A) shows the perturbation vector field associated with the matrix perturbation $\mathbf{P}$ that maximizes the expectation value of $\Psi$. Such a vector field is superimposed on the invariant measure of the unperturbed system. Panel (B) shows the linear response with respect to the perturbation $\mathbf{P}$. Panel (C) is the difference between the diagonal diffusion tensors obtained from $\mathbf{M}+\varepsilon \mathbf{P}$ and $\mathbf{M}$.}
\end{figure}

In Fig.~\ref{fig:dwp2_resp_lowres}, we redo the lower-resolution analyses for obtaining a perturbation matrix $\mathbf{P}$ such that the expectation value of the observable $\Psi$ defined in Eq.~\eqref{eq:observable} is maximized. Panel (A) shows the linear response to the perturbation $\mathbf{P}$. Panel (B) is the vector field associated with $\mathbf{P}$. Panel (C) is the diffusion introduced by $\mathbf{P}$. We observe that the results compared to the higher resolution plots in Fig.~\ref{fig:dwp1_depth14} are stable.

\subsection{Dependence on perturbation strength $\varepsilon$}\label{app:dependence_epsilon}

\begin{figure}[H]
\centering
\includegraphics[scale=0.29]{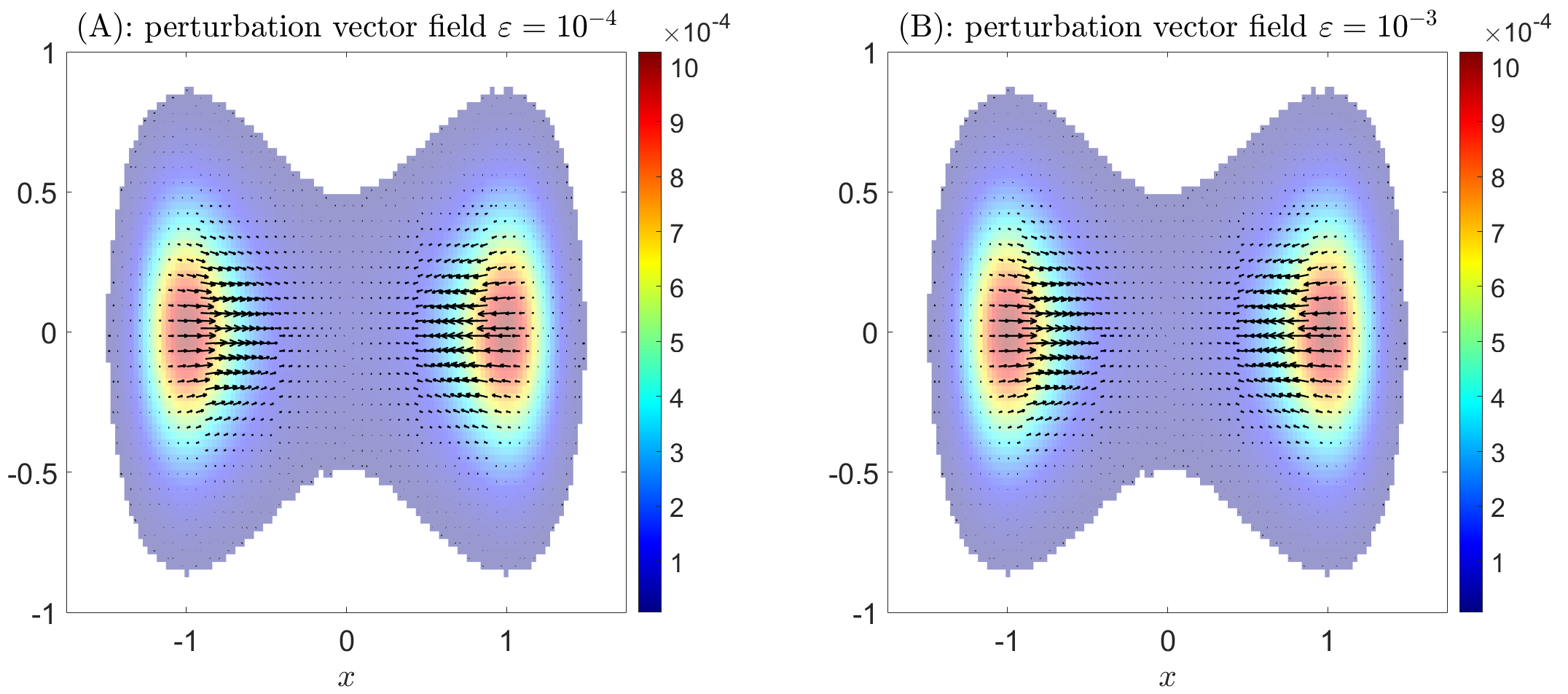}
\caption{\label{fig:epsilon_Dependence}\textbf{Perturbation reconstruction: $\varepsilon$ dependence.} We show the reconstruction of the perturbation vector field associated with $\mathbf{P}$, where $\mathbf{P}$ maximizes the expectation value of the observable $\Psi$ defined in Eq.~\eqref{eq:observable}. Panel (A) corresponds to the reconstruction algorithm using $\varepsilon = 10^{-4}$ and Panel (B), $\varepsilon = 10^{-3}$.}
\end{figure}

The results of the perturbation matrices and their vector field reconstruction do not depend on $\varepsilon$. Here we show the vector field reconstruction of the matrix $\mathbf{P}$ that was obtained to optimize the expectation value of the observable $\Psi$ defined in Eq.~\eqref{eq:observable}. Two other values of $\varepsilon$ are considered on top of $\varepsilon=0.01$ shown in the main text. Panel (A) in Fig.~\ref{fig:epsilon_Dependence} shows the reconstruction for $\varepsilon = 10^{-3}$ and Panel (B) shows the same for $\varepsilon=10^{-4}$. We observe that the differences are minimal.

\section{Undoing the rotation}\label{app:additive_reversibilization_check}

The target of this appendix is two-fold. First, we want to test the finite-difference formula shown in Eq.~\eqref{eq:drift_recons}. Second, we want to find out what vector field is associated with the perturbation matrix $\mathbf{P}_r$ defined in Eq.~\eqref{eq:reversibilized_perturbation}. We recall that that such matrix contains the negative of the anti-symmetric part of the matrix $\mathbf{M}$, where symmetry is understood in the inner product weighted by the entry-wise inverse of the invariant vector of $\mathbf{M}$. One can, thus, speculate that the vector field associated with $\mathbf{P}_r$ is the reverse of the applied rotation $\mathbf{R}$ in Eq.~\eqref{eq:rotation}, in the context of the two-dimensional double well potential.

\begin{figure}[H]
\centering
\includegraphics[scale=0.29]{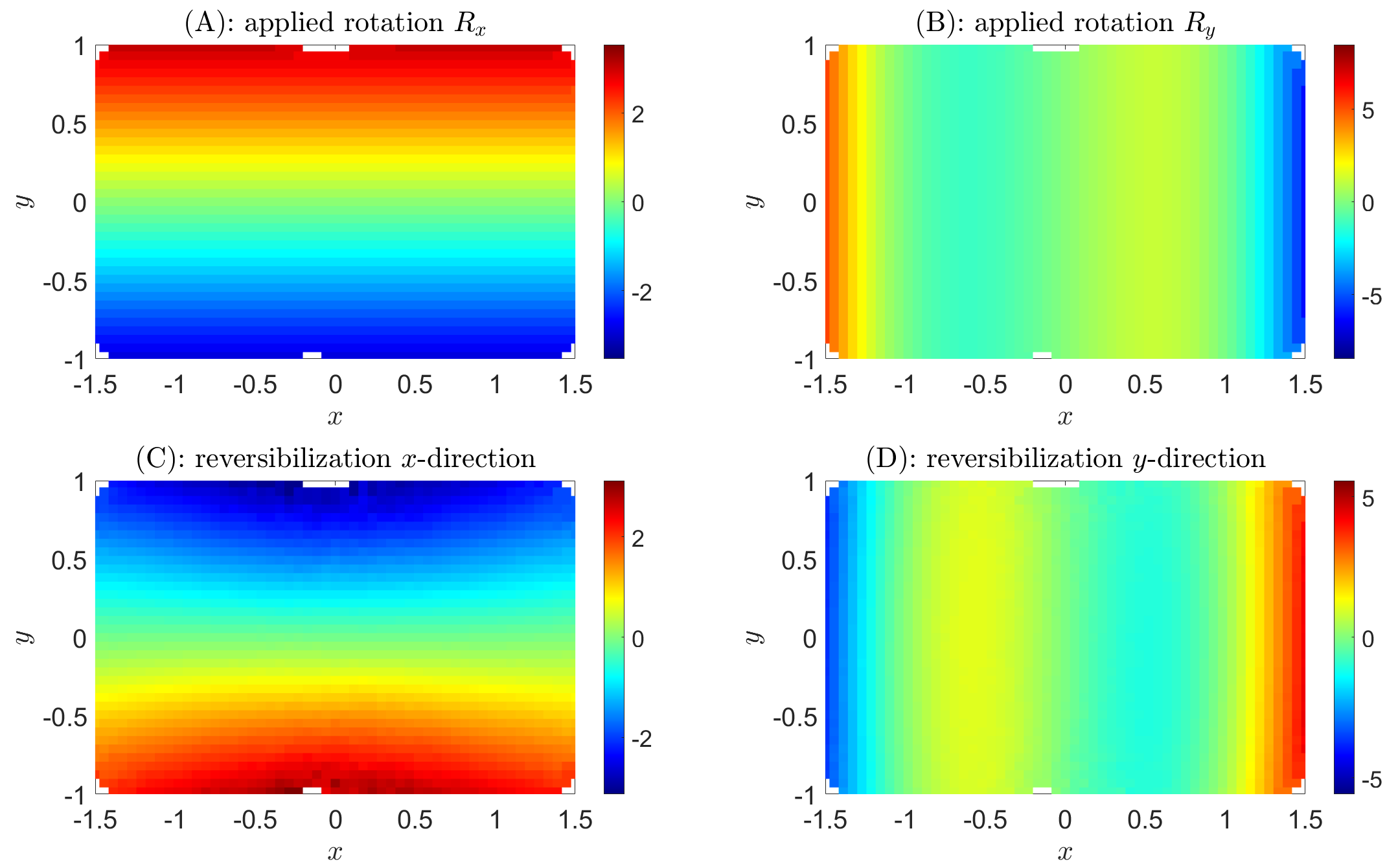}
\caption{\label{fig:comparison_rev}\textbf{Reconstructing the rotation.} Panels (A) and (B) show the $x$ and $y$ components of the applied rotation $\mathbf{R}$ in Eq.~\eqref{eq:rotation}. Panels (C) and (D) shows the $x$ and $y$ components of formula Eq.~\eqref{eq:drift_recons} applied to $\mathbf{P}_r$ defined in Eq.~\eqref{eq:reversibilized_perturbation}.}
\end{figure}

In Fig.~\ref{fig:comparison_rev}(A) and (B) we show the $x$ and $y$ components of the vector field $\mathbf{R}$ defined in Eq.~\eqref{eq:rotation}. We then applied the formula in Eq.~\eqref{eq:drift_recons} to reconstruct the vector field associated with $\mathbf{P}_r$. Panels (C) and (D) shows the $x$ and $y$ components respectively. Indeed, a qualitative--- and also quantitative, see colorbars--- look to the figure shows that the vector field does as posited in the previous paragraph. This also shows the effectiveness of formula Eq.~\eqref{eq:drift_recons}.

\section{Projection method for constraints}\label{sec:projection_space}

In this appendix, we show that the formula in Eq.~\eqref{eq:lagrange_matrix_formula} is equivalent to a projection onto the null-space of a suitably constructed linear operator. We consider the problem of optimizing the linear functional of the form:
\begin{equation}\label{eq:linear_functional}
f(\PP) =  \sum_{i=1}^N\sum_{j=1}^N C_{ij} P_{ij} = \widehat{\mathbf{C}}\cdot \widehat{\mathbf{P}}
\end{equation}
where $\widehat{\phantom{\mathbf{Aa}}}$ denotes the vectorization of a matrix, $\mathbf{C}$ is a given matrix and we impose the following constraints:
\begin{align*}
\text{(C1)} &\quad \|\PP\|_F = 1 \\
\text{(C2)} &\quad \mathbf{1}^\top \PP = 0^\top \quad \text{(columns of } \PP \text{ sum zero)} \\
\text{(C3)} &\quad P_{ij} = 0 \quad \text{( if } M_{ij} = 0 \text{)}\\
\text{(C4)} &\quad \PP \uu = \mathbf{0} \quad \text{(columns lie in the nullspace of } \uu^\top \text{)}
\end{align*}
In \cite{Antown2018} they prove the existence of an optimizer to constraints (C1)-(C2)-(C3) by reformulating the problem in vectorized form:
\begin{lemma}[Antown et al. 2018]
    The constraint (C2) can be written in the form of $\AAA \widehat{\PP} = \mathbf{0}$, where $\AAA$ is a $N\times N^2$ matrix and $\widehat{\PP}$ is the vectorization of the matrix $\PP$ and $\AAA$ is defined as:
    \begin{equation}
       \AAA =  \mathbf{I}_N \otimes \mathbf{1}^{\top},
    \end{equation}
    where $\mathbf{I}_N$ is the $N$-dimensional identity matrix.
\end{lemma}

We show an analogous formulation for (C4):
\begin{lemma}
    The constraint (C4) can be written in the form of $\BB \widehat{\PP} = \mathbf{0}$, where $\BB$ is a $N\times N^2$ matrix defined as:
    \begin{equation}
       \BB =  \left[\mathbf{I}_N \otimes \mathbf{u}^{\top}\right]\mathbf{K},
    \end{equation}
    where $\mathbf{K}$ is the $N^2\times N^2$ matrix defined as:
    \begin{equation}
        K_{(i-1)N+j,(j-1)N+i} = 1,
    \end{equation}
    for $i,j=1,\ldots,N$. This matrix satisfies $\mathbf{K}\widehat{\PP} = \widehat{\PP^{\top}}$.
\end{lemma}
\begin{proof}
    Since $\PP\uu = \mathbf{0}$ if and only if $\uu^{\top}\PP^{\top}=\mathbf{0}^{\top}$, if and only if, $\PP\uu = \widehat{\uu^{\top}\PP^{\top}}$ and:
    \begin{equation}
        \PP \uu = \widehat{\uu^{\top}\PP^{\top}} = \widehat{\uu^{\top}\PP^{\top}\mathbf{I}_N} = \left[ \mathbf{I}_N \otimes \uu^{\top} \right]\widehat{\PP^{\top}} = \left[ \mathbf{I}_N \otimes \uu^{\top} \right]\mathbf{K}\widehat{\PP}.
    \end{equation}
    
We are left with showing that $\mathbf{K}\widehat{\PP} = \widehat{\PP^{\top}}$. Let $\mathbf{v} = \widehat{\PP}$. Note that for every $1\leq i,j \leq N$, we have that $P_{ij} = v_{(j-1)N+i}$. We need to prove that this indexation is well defined. For this it is enough to show that the mapping $\iota: \{1,\ldots,N\}\times \{1,\ldots,N\} \rightarrow \{1,\ldots,N^2\}$ defined by $\iota(i,j)=(j-1)N+i$ is invertible. To prove the one-to-one condition, let us assume that $\iota(i_1,j_1)=\iota (i_2,i_2)$, for $i_1,i_2,j_1,j_2$ in $\{1,\ldots,N\}$. We then have:
\begin{equation}
(j_1-1)N+i_1 = (j_2-1)N+i_2 \Leftrightarrow (j_1 - j_2)N = i_2 - i_1,    
\end{equation}
from where $i_2-i_1$ must be divisible by $N$. However, since $i_1,i_2$ are positive integers smaller than $N$, this can only happen if $i_1=i_2$. This implies that $j_1=j_2$.

To prove surjectivity, let $k$ be in $\{1,\ldots, N^2\}$ and define $j = \lfloor(k-1)/N\rfloor+1$ and $i=k-(j-1)N$. Then,
\begin{align}
    \iota(i,j) &= \left(\left(\left\lfloor\frac{k-1}{N}\right\rfloor+1\right)  -1\right)N + k-\left(\left(\left\lfloor\frac{k-1}{N}\right\rfloor+1\right)-1\right)N = k.
\end{align}
Furthermore, $1\leq i,j\leq N$ for any $k$.

Let $\mathbf{w}=\mathbf{K}\widehat{\PP}$. We have to show that $w_{(j-1)N+i} = P_{ji}$. By the definition of $\mathbf{K}$ we have:
\begin{equation}
    w_{(j-1)N+i} = \sum_{k=1}^N K_{(j-1)N+i,k}v_k = v_{(i-1)N+j} = P_{ji},
\end{equation}
and this equality holds for every $1\leq i,j\leq N$. This completes the proof.

\end{proof}
It remains to prove that the intersection of the nullspaces of $\mathbf{A}$ and $\mathbf{B}$ is nontrivial:
\begin{lemma}
Let $N \geq 3$, then the set $\mathcal{V} = \left\{ \PP \in \mathbb{R}^{N\times N} : \mathbf{1}^\top \PP = \uu^{\top}\PP^{\top} = \mathbf{0}^{\top} \right\}$ is a vector space of dimension at least $N^2 - 2N$.
\end{lemma}
\begin{proof}
    Let us block-wise define the matrix $\mathbf{S}$ in $\mathbb{R}^{2N\times N^2}$ as:
    \begin{equation}\label{eq:matrix_S}
        \mathbf{S} = \begin{bmatrix} \mathbf{A} \\ \mathbf{B} \end{bmatrix}.
    \end{equation}
    Then, $\ker(\mathbf{S})  = \ker(\AAA)\cap \ker(\BB)$. Indeed, $\xx$ is in $\ker(\mathbf{S})$ if and only if $\mathbf{A}\xx = \mathbf{0}$ and $\mathbf{B}\xx = \mathbf{0}$, so that $\xx$ is simultaneously in $\ker(\mathbf{A})$ and $\ker(\mathbf{B})$. Finally, because $\mathbf{S}$ has an image space within $\mathbb{R}^{2N}$ and an input space within $\mathbb{R}^{N^2}$, we have that $\dim \left( \ker(\mathbf{S}) \right) = N^2  - \mathrm{rank}(\mathbf{S})$. But $\mathrm{rank}(\mathbf{S})\leq 2N$. 
    Hence, $\mathcal{V}$ is a vector space with dimension at least $N^2  - 2N$.
\end{proof}

At this stage, the kernel of the matrix $\mathbf{S}$ defined in Eq.~\eqref{eq:matrix_S} serves for the projection of optimization algorithms for any Markov matrix with $M_{ij}>0$ for all $i,j=1,\ldots,N$. But not only for linear functionals like that in Eq.~\eqref{eq:linear_functional}, but also for nonlinear ones. 

To obtain the exact optimizing matrix $\mathbf{P}$ in Eq.~\eqref{eq:linear_functional}, we are left with constraining the solutions so that the sparsity of $\mathbf{M}$ is respected, namely, constraint (C3). Let us define $I = \{ i_j \}_{j=1}^{\# I} =  \left\{ k \in \{1,\ldots,N \} : \left[\widehat{\mathbf{M}}\right]_{k} > 0 \right\}$, where $\#I$ is the number of non-zero entries in $\mathbf{M}$. We also introduce the reduced vectors $\mathbf{c}$ and $\mathbf{p}$  and reduced matrices $\AAA^r$, $\BB^r$ and $\mathbf{S}^r$ defined by:
\begin{subequations}
    \begin{align}
        c_j &= \left[ \widehat{\CC} \right]_{i_j} , \label{eq:little_c}\\
        p_j &= \left[ \widehat{\PP} \right]_{i_j} , \label{eq:little_p}\\
        A^r_{kj} &= A_{kj} ,\\
        B^r_{kj} &= B_{kj}, \\
        \mathbf{S}^r &= \begin{bmatrix}
            \mathbf{A}^r \\ \mathbf{B}^r
        \end{bmatrix}, \label{eq:reduced_S}
    \end{align}
\end{subequations}
where $i_j$ is in $I$ and $j=1,\ldots,\# I$ and $k=1\ldots, N$. Notice that $\mathbf{c}$ and $\mathbf{p}$ are in $\mathbb{R}^{\# I}$, $\AAA^r$ and $\BB^r$ are in $\mathbb{R}^{N\times \# I}$ and $\mathbf{S}^r$ in $\mathbb{R}^{2N\times \# I}$. Then we have the following equality:
\begin{equation}
    f(\PP) = \widehat{\CC}\cdot \widehat{\PP} = \mathbf{c}\cdot \mathbf{p},
\end{equation}
hence we have to optimize  $\mathbf{p}$ over the kernel of $\mathbf{S}^r$. We have, however,
\begin{equation}
    \dim \left( \ker(\mathbf{S}^r) \right) = \#I  - \mathrm{rank}(\mathbf{S}^r).
\end{equation}

% \begin{lemma}\label{lemma4}
%     Let $\mathbf{M}$ be a mixing $N\times N$ Markov matrix with a stationary vector $\uu = (u_i)$ such that $u_i>0$ for every $i=1,\ldots,N$. Then, $\#I\geq2N$.
% \end{lemma}
% \begin{proof}
% First of all, all the columns of $\mathbf{M}$ must add to one and $M_{ij}>0$ for every $i,j=1,\ldots,N$, so each column has at least one non-zero entry. Thus $N\leq \#I$. Secondly, if column $j_0$ has a single non-zero entry $M_{i_0,j_0}>0$, this implies that $M_{i_0,j_0} = 1$. Indeed,
% \begin{equation}
%     1 = \sum_{i=1}^NM_{ij_0} = M_{i_0j_0}.
% \end{equation}
% Hence, each column of $\mathbf{M}$ has at least two nonzero entries, meaning that $\#I\geq 2N$.
% \end{proof}

So we have to show that $\#I>2N\geq\mathrm{rank}(\mathbf{S}^r) $ to guarantee a nontrivial projection space. Since we require that $M_{ij}<1$ for every $i,j=1,\ldots,N$, each column must have at least two non-zero entries, which implies that $\#I\geq 2N$. Hence, it is possible that the kernel of $\mathbf{S}^r$ has zero dimension. In any case, provided that such a kernel is nontrivial, the solutions of the optimization problem of Eq.~\eqref{eq:linear_functional} are given by the orthogonal projection of $\mathbf{c}$ in Eq.~\eqref{eq:little_c} onto the kernel of $\mathbf{S}^r$ in Eq.~\eqref{eq:reduced_S}. That would provide the vector $\mathbf{p}$ in Eq.~\eqref{eq:little_p} which, upon reordering, gives the desired $\mathbf{P}$.

\subsection{Comparison between the two methods}

In the main text we provided the solution to the Lagrange multiplier equation for the optimization of Eq.~\eqref{eq:linear_functional}, shown in formula \eqref{eq:antown_expect}. We call this Method~1. In this appendix, we also show that the solutions to this optimization problem can be viewed as the orthogonal projection of the data matrix onto the kernel of a suitable functional--- described in details in the previous section---. We shall call this Method~2. In principle, both methods should provide the same results, but the radically different implementation and conceptual basis is prone to numerical disagreements. Here we show that both methods give the same results. Method~1 and Method~2 are applied to the entropy production minimization problem of Eq.~\eqref{eq:min_ent_prod}. This, in addition, confirms that the methods indeed minimize the entropy production. 

We randomly initialized $K = 10^4$ matrices of size $50\times 50$ with a diagonal dominance of $100\%$--- the diagonal element is the sum of the non-diagonal entries in its column--- and a sparsity of $50\%$. The entries are randomly sampled from a uniform distribution in $[0,1]$ and normalized to obtain Markov matrices. Non-mixing matrices--- for instance, if a diagonal element was equal to unity--- are discarded and randomly replacing them so that the ensemble consists of $K=10^4$ members. The factors of diagonal dominance and sparsity control the spectral gap--- defined as $1-|\lambda_1|$, where $\lambda_1$ is the subdominant eigenvalue---. We chose the values so that they resemble the sparsity and diagonal dominance of the matrices considered earlier in this study. A mean spectral gap of $0.31\pm 0.03$ was obtained.

\begin{figure}[H]
\centering
\includegraphics[scale=0.29]{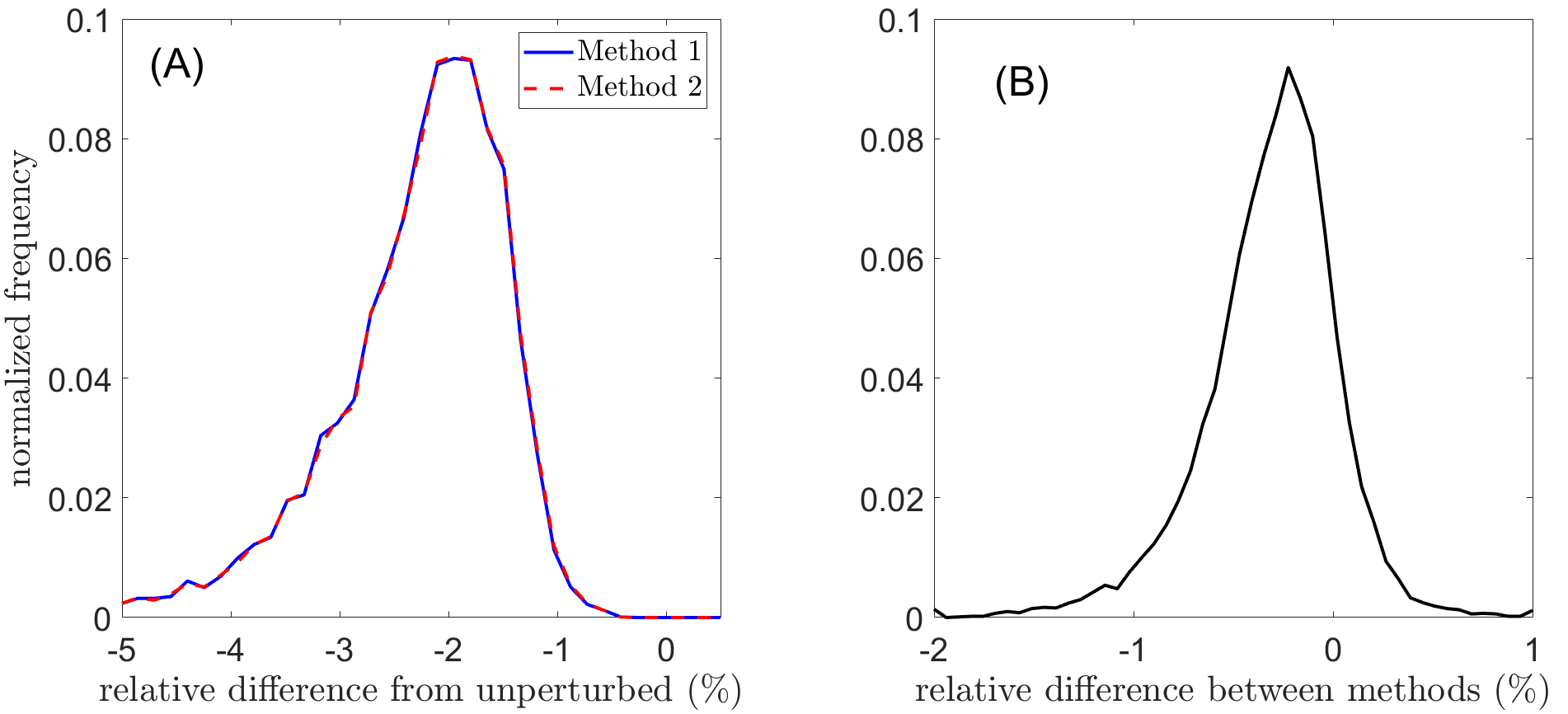}
\caption{\label{fig:comparison}\textbf{Comparison between optimization methods.} Panel (A): normalized histograms for the relative difference between $s\left( \mathbf{M}^{(k)} + \varepsilon\mathbf{P}_{i}^{(k)}\right)$ and $s\left( \mathbf{M}^{(k)} \right)$ for both methods shown in the legend. Panel (B): normalized histograms of the relative difference between methods in the minimization of entropy production.}
\end{figure}

Then Method~1 and Method~2 are applied to the matrices to obtain perturbation matrices $\mathbf{P}_i^{(k)}$, for $k=1,\ldots,K$ and $i=1,2$--- the subscript $i$ refers to the method used---. Then, for a value of $\varepsilon = 10^{-3}$, we computed the entropy production $s\left( \mathbf{M}^{(k)} + \varepsilon\mathbf{P}_{i}^{(k)}\right)$ and computed the relative difference with respect to $s\left( \mathbf{M}^{(k)} \right)$ for each random initialization and method. Then we computed the (normalized) histogram shown in Fig.~\ref{fig:comparison}(A). Indeed, this shows that both methods actually minimize the entropy production consistently. Secondly, we computed the histogram of relative difference between methods for each random initialization of the matrix $\mathbf{M}^{(k)}$. The normalized histogram in Fig.~\ref{fig:comparison}(B) shows that both methods provide different results, but that their difference is small and here we attribute it to numerical issues.

\paragraph{Data availability}\quad The data used in this work is available upon request.

\paragraph{Conflict of interests}\quad The authors declare that they have no competing interests.
 
\AtNextBibliography{\footnotesize}
\printbibliography

\end{document}